\documentclass[12pt]{article}
\usepackage{amsfonts,amsmath}
\def \beq {\begin{eqnarray}}
\def \eeq {\end{eqnarray}}
\def \beqn {\begin{eqnarray*}}
\def \eeqn {\end{eqnarray*}}

\newcommand{\halmos}{\rule{1ex}{1.4ex}}

\newcounter{for}[section]

\newtheorem{itlemma}{Lemma}[section]
\newtheorem{itproposition}[itlemma]{Proposition}
\newtheorem{theorem}[itlemma]{Theorem}
\newtheorem{itcorollary}[itlemma]{Corollary}
\newtheorem{itremark}[itlemma]{Remark}
\newtheorem{itremarks}[itlemma]{Remarks}
\newtheorem{itdefinition}[itlemma]{Definition}
\newtheorem{itexample}[itlemma]{Example}

\newenvironment{fact}{\begin{itfact}\rm}{\end{itfact}}
\newenvironment{claim}{\begin{itclaim}\rm}{\end{itclaim}}
\newenvironment{lemma}{\begin{itlemma}}{\end{itlemma}}
\newenvironment{remark}{\begin{itremark}\it}{\end{itremark}}
\newenvironment{remarks}{\begin{itremarks} \rm}{\end{itremarks}}
\newenvironment{corollary}{\begin{itcorollary}}{\end{itcorollary}}
\newenvironment{proposition}{\begin{itproposition}}{\end{itproposition}}
\newenvironment{definition}{\begin{itdefinition}\rm}{\end{itdefinition}}
\newenvironment{example}{\begin{itexample}\rm}{\end{itexample}}
\newenvironment{proof}{\noindent {\em Proof}.\ \
}{\hspace*{\fill}$\halmos$\medskip}
\newcommand{\be}[1]{\addtocounter{for}{1} \begin{equation}\label{#1}}
\newcommand{\ee}{\end{equation}}
\newcommand{\bl}[1]{\begin{lemma}\label{#1}}
\newcommand{\br}[1]{\begin{remark}\label{#1}}
\newcommand{\brs}[1]{\begin{remarks}\label{#1}}
\newcommand{\bt}[1]{\begin{theorem}\label{#1}}
\newcommand{\bd}[1]{\begin{definition}\label{#1}}
\newcommand{\bp}[1]{\begin{proposition}\label{#1}}
\newcommand{\bc}[1]{\begin{corollary}\label{#1}}
\newcommand{\bfact}[1]{\begin{fact}\label{#1}}
\newcommand{\bex}[1]{\begin{example}\label{#1}}
\newcommand{\ec}{\end{corollary}}
\newcommand{\efact}{\end{fact}}
\newcommand{\eex}{\end{example}}
\newcommand{\el}{\end{lemma}}
\newcommand{\er}{\end{remark}}
\newcommand{\ers}{\end{remarks}}
\newcommand{\et}{\end{theorem}}
\newcommand{\ed}{\end{definition}}
\newcommand{\ep}{\end{proposition}}
\newcommand{\epr}{\end{proof}}
\newcommand{\bpr}{\begin{proof}}
\newcommand{\bcl}[1]{\begin{claim}\label{#1}}
\newcommand{\ecl}{\end{claim}}
\newcommand{\ecs}{\end{corollary}}
\newcommand{\eers}{\end{exercise}}
\newcommand{\eexs}{\end{example}}
\newcommand{\eems}{\end{example}}
\newcommand{\els}{\end{lemma}}
\newcommand{\eles}{\end{lemmaex}}
\newcommand{\ets}{\end{theorem}}
\newcommand{\eds}{\end{definition}}
\newcommand{\eps}{\end{proposition}}

\newcommand{\bi}{\begin{itemize}}
\newcommand{\ei}{\end{itemize}}
\newcommand{\ben}{\begin{enumerate}}
\newcommand{\een}{\end{enumerate}}

\def\vbar{\mathchoice{\vrule height6.3ptdepth-.5ptwidth.8pt\kern-.8pt}
   {\vrule height6.3ptdepth-.5ptwidth.8pt\kern-.8pt}
   {\vrule height4.1ptdepth-.35ptwidth.6pt\kern-.6pt}
   {\vrule height3.1ptdepth-.25ptwidth.5pt\kern-.5pt}}
\def\fudge{\mathchoice{}{}{\mkern.5mu}{\mkern.8mu}}
\def\bbc#1#2{{\rm \mkern#2mu\vbar\mkern-#2mu#1}}
\def\bbb#1{{\rm I\mkern-3.5mu #1}}
\def\bba#1#2{{\rm #1\mkern-#2mu\fudge #1}}
\def\bb#1{{\count4=`#1 \advance\count4by-64 \ifcase\count4\or\bba A{11.5}\or
   \bbb B\or\bbc C{5}\or\bbb D\or\bbb E\or\bbb F \or\bbc G{5}\or\bbb H\or
   \bbb I\or\bbc J{3}\or\bbb K\or\bbb L \or\bbb M\or\bbb N\or\bbc O{5} \or
   \bbb P\or\bbc Q{5}\or\bbb R\or\bbc S{4.2}\or\bba T{10.5}\or\bbc U{5}\or
   \bba V{12}\or\bba W{16.5}\or\bba X{11}\or\bba Y{11.7}\or\bba Z{7.5}\fi}}

\def \qed {{\hspace*{\fill}$\halmos$\medskip}}

\def \L {{\cal{L}}}

\def \A {{\cal{A}}}

\def \G {{\cal{G}}}
\def \C {{\cal{C}}}
\def \D {{\cal{D}}}

\newcommand{\ba}[1]{\addtocounter{for}{1} \begin{eqnarray}\label{#1}}
\newcommand{\ea}{\end{eqnarray}}

\def\sqr#1#2{{\vcenter{\vbox{\hrule height .#2pt
                             \hbox{\vrule width .#2pt height#1pt \kern#1pt
                                   \vrule width .#2pt}
                             \hrule height .#2pt}}}}

\def\pmb#1{\setbox0=\hbox{#1}%
   \kern-.025em\copy0\kern-\wd0
   \kern.05em\copy0\kern-\wd0
 \kern-.025em\raise.0433em\box0 }
\def\sqr#1#2{{\vcenter{\vbox{\hrule height.#2pt
     \hbox{\vrule width.#2pt height#1pt \kern#1pt
   \vrule width.#2pt}\hrule height.#2pt}}}}

\def\B{{\mathcal B}}

\def \I{{\mathcal I}}
\def\N{{\mathbb N}}
\def\Z{{\mathbb Z}}
\def\R{{\mathbb R}}

\def\var{\text{var}}

\def\bs{\backslash}
\def\reff#1{(\ref{#1})}

\def \ind {\hbox{1\hskip -3pt I}}
\newcommand {\cro}[1] {\left[ {#1} \right]}
\newcommand {\acc}[1] {\left\{ {#1} \right\}}
\newcommand {\pare}[1] {\left( {#1} \right)}

\begin{document}

\title{Shape transition under excess self-intersections for transient random walk.}
\author{Amine Asselah \\ Universit\'e Paris-Est\\
Laboratoire d'Analyse et de Math\'ematiques Appliqu\'ees\\
UMR CNRS 8050\\ amine.asselah@univ-paris12.fr}
\date{}
\maketitle
\begin{abstract}
We reveal a shape transition for
a transient simple random walk forced to realize 
an excess $q$-norm of the local times,
as the parameter $q$ crosses the value $q_c(d)=\frac{d}{d-2}$. 
Also, as an application of our approach,
we establish a central limit theorem
for the $q$-norm of the local times in dimension 4 or more.
\\ \\
{\em Abstract in French:}
Nous d\'ecrivons un ph\'enom\`ene de transition de forme
d'une marche al\'eatoire transiente forc\'ee \`a r\'ealiser
une grande valeur de la norme-$q$ du temps local, 
lorsque le param\`etre $q$ 
traverse la valeur critique $q_c(d)=\frac{d}{d-2}$.
Comme application de notre approche, nous \'etablissons un th\'eor\`eme
de la limite centrale 
pour la norme-$q$ du temps local en dimension 4 et plus.

\end{abstract}

{\em Keywords and phrases}: self-intersection local times, large deviations, random walk.

{\em AMS 2000 subject classification numbers}: 60K35, 82C22,
60J25.

{\em Running head}: Shape transition under excess self-intersection.

\section{Introduction} \label{sec-intro}
We consider a simple random walk
$\{S(n),\ n\in \N\}$ on $\Z^d$, starting at the origin. 
For any set $A$, we denote
by $\ind_A$ the indicator of $A$, and consider the local times of the walk 
$\{l_n(z),\ z\in \Z^d\}$ given by
\be{intro.2}
l_n(z)=\ind_{\acc{S(0)=z}}+\dots+\ind_{\acc{S(n-1)=z}}.
\ee
For a real $q>1$, we form the sum of the $q$-th power of
the local times
\be{intro.1}
||l_n||_q^q=\sum_{z\in \Z^d} l_n(z)^q.
\ee
When $q$ is integer, $||l_n||_q^q$ can be written in
terms of the $q$-fold self-intersection local times of a random walk.
For instance, when $q=2$
\[
||l_n||_2^2=n+2\sum_{0\le i<j<n} \ind_{\acc{S(i)=S(j)}}.
\]
For $q$ positive real, 
we still call $||l_n||_q^q$ the $q$-fold self-intersection
local times.

In dimension three and more, Becker and K\"onig \cite{BK} have shown that
there are positive constants, say $\kappa(q,d)$, such that almost surely
\be{intro.3}
\lim_{n\to\infty} \frac{||l_n||_q^q}{n}=\kappa(q,d).
\ee
Here, we are concerned with estimating the deviations
of $||l_n||_q^q$ away from its mean. That is, if $P$ denotes
the law of the walk started at 0, we give estimates for
\be{question}
P\pare{||l_n||_q^q-E[||l_n||_q^q]\ge \xi n}.
\ee
for $\xi$ positive, and $n$ going to infinity. 

There is a rich literature concerning the two-fold self-intersection
local times. The reason being that $||l_n||_2$ is a natural object in
quantum-field theory (see \cite{aizenman}, \cite{felder-frolich}
and \cite{varadhan} for instance),
as well as in the statistical physics of polymers (see 
\cite{edwards}, \cite{brydges-slade}
and \cite{BS} for instance). However $||l_n||_q$ for 
$q\in \R\bs \N$ has no such direct links with physics. It comes up
naturally in studying large and moderate deviations for
{\it random walk in random sceneries} (see \cite{AC05} and
\cite{FMV}).

Now, in the large deviations results for the two-fold self-intersection
of a transient random walk (see \cite{AC05,A06,Chen08,A07})
two strategies have a distinguished r\^ole. 
\begin{itemize}
\item Strategy A: the walk visits of the order of
$(\xi n)^{1/q}$-times, finitely
many sites in a ball of bounded radius. For a transient walk,
the number of visits
of a bounded domain is bounded by a geometric variable. Thus, strategy A
costs of the order of $\exp(-O((n\xi)^{1/d}))$, where we use the notation
$y_n=O(x_{n})$ for two positive sequences $\{x_{n},y_n,\ n\in \N\}$,
to mean that there is $K>0$ such that $0\le y_n\le K x_{n}$.
\item Strategy B: the walk visits of the order of $\xi^{1/(q-1)}$-times,
about $n/\xi^{1/(q-1)}$ sites. Presumably, the walk stays, a time $n$,
in a ball
of volume $n/\xi^{1/(q-1)}$. The cost of staying a time $n$ within
a ball of radius $r_n\ll \sqrt n$ is about $\exp(-O(n/r_n^2))$, so that
strategy B costs of the order of $\exp(-O(n^{1-\frac{2}{d}}\xi^{
\frac{2}{d(q-1)}}))$.
\end{itemize}
When $q=2$, \cite{AC05,A07} have shown that strategy A is adopted in
$d\ge 5$, whereas \cite{A06} (see also Chapter 8.4 of \cite{Chen08})
suggests that strategy B is adopted in $d=3$. 

To summarize in words our main finding, assume
$d\ge 3$, fix $\xi>0$ and look at typical paths
realizing $\{||l_n||_q^q-E[||l_n||_q^q]\ge \xi n\}$. 
As we increase $q$, we step on 
a value, $q_c(d)$, above which our large deviation event is 
realized by strategy A, and below which it is realized by
strategy B. The critical value $q_c(d)=\frac{d}{d-2}$ is 
obtained as we equal the costs of strategies A and B.

Note that $q_c(d)$ is a well known number: if $q$ is integer,
then $q$ independent simple random walks, on $\Z^d$, intersect
infinitely often if and only if $q<q_c(d)$ (see for instance
\cite{legall86} Proposition 7.1, and \cite{lawler} Section 4.1).

Let us now describe, in mathematical terms, this shape transition. 
The first theorem deals with the {\it sub-critical regime} $q<q_c(d)$.

\bt{theo-d3} Assume dimension $d\ge 3$. Then, for $q$ and $d$
such that $1<q<\frac{d}{d-2}$, there are 
constants $c_1^\pm(q,d)>0$ such that for $\xi\ge 1$, and $n$ large enough
\be{intro.4}
\exp\pare{-c_1^-(q,d) \xi^{\frac{2}{d}(\frac{1}{q-1})}n^{1-\frac{2}{d}}}
\le\ P\pare{ ||l_n||_q^q-E[||l_n||_q^q]\ge \xi n}\ \ \le
\exp\pare{-c_1^+(q,d) \xi^{\frac{2}{d}(\frac{1}{q-1})}n^{1-\frac{2}{d}}}.
\ee
Moreover, in this regime the sites visited more than some large constant do not
contribute to realizing the excess self-intersection. In other words,
\be{intro.5}
\limsup_{A\to\infty}\limsup_{n\to\infty} \frac{1}{n^{1-\frac{2}{d}}}
\log\ P\pare{\sum_{z\in \Z^d}\ind_{\acc{l_n(z)>A}}l_n(z)^q\ge
\xi n}=-\infty.
\ee
\et
Our second theorem deals with the {\it super-critical regime} $q>q_c(d)$.
\bt{theo-d4} Assume dimension $d\ge 3$. For $q$ and $d$ such that
$q>\frac{d}{d-2}$, there are constants $c_2^\pm(q,d)>0$ 
such that for $\xi\ge 1$, and $n$ large enough 
\be{intro.6}
\exp\pare{-c_2^-(q,d) (\xi n)^{1/q}}\le\ \ 
P\pare{ ||l_n||_q^q-E[||l_n||_q^q]\ge \xi n}\ \ \le \exp\pare{-c_2^+(q,d)
(\xi n)^{1/q}}.
\ee
Moreover, the sites visited much less than $n^{\frac{1}{q}}$ do not
contribute to realizing the excess self-intersection. In other words,
\be{intro.7}
\limsup_{\epsilon\to0}\limsup_{n\to\infty} \frac{1}{n^{1/q}}
\log\ P\pare{\sum_{z\in \Z^d}\ind_{\acc{l_n(z)<\epsilon n^{1/q}}}
l_n(z)^q\ge E[||l_n||_q^q]+\xi n}=-\infty.
\ee
\et
\br{rem-growth}
In Theorems~\ref{theo-d3} and \ref{theo-d4},
we could take $\xi$ to grow with $n$. 
The only (necessary) bound on $\xi_n$ comes from
the bound $||l_n||_q\le n$ which imposes that $\xi_n\le n^{q-1}$.
The proofs are written with general $\xi_n\ge 1$.
\er

The next result deals with the contribution of some level
sets of the local times to deviation on a much larger
scale than the mean, and can be obtained by the same approach
yielding Theorem~\ref{theo-d4}.
We include it in this form since it can be of independent interest,
while showing the possibilities offered by our approach.
Also, it generalizes Lemma 1.8 of \cite{AC05}.
\bl{lem-level}
Assume $d\ge 3$ and $q\ge q_c(d)$. Choose $a,b>0$ such that
$1<a<1+b(q-1)$. Then, for any $\epsilon>0$, and $n$ large enough
\be{level-ineq1}
P\pare {\sum_{z\in \Z^d}\ind_{\acc{l_n(z)<n^b}}l_n(z)^q\ge n^a}\le
e^{-n^{\zeta(q,a,b)-\epsilon}}\text{  with  }
\zeta(q,a,b)=b+\frac{1}{q_c(d)}(a-qb).
\ee
\el

\br{rem-xin}
Our approach is not suited to studying
small $\xi_n$ for reasons explained later in Remark~\ref{rem-realq}.
However, when $1>\xi_n\ge n^{-\delta}$, for some positive
$\delta$ small enough, our approach yields
a constant $c_1$ such that for $q<q_c(d)$
\be{eq-xin1}
P\pare{ ||l_n||_q^q-E[||l_n||_q^q]\ge \xi_n n}\ \ \le 
\exp\pare{- c_1\xi_n^{\frac{2}{d}(\frac{q}{q-1})}n^{1-\frac{2}{d}}}.
\ee
When $q>q_c(d)$, we have a constant $c_2$ such that
\be{eq-xin2}
P\pare{ ||l_n||_q^q-E[||l_n||_q^q]\ge \xi_n n}\ \ \le 
\exp\pare{- c_2\xi_n^{\frac{1}{q}+\frac{2}{d}}n^{\frac{1}{q}}}.
\ee
We believe that the powers of $\xi_n$ in
\reff{eq-xin1} and \reff{eq-xin2} are not optimal. 
However, \reff{eq-xin1} and 
\reff{eq-xin2} are useful in deriving a central limit
theorem stated in Theorem~\ref{theo-clt}.
\er

Our initial goal was to improve the main result of~\cite{A06}, which
states that in dimension 3, 
there is $\underline\chi>0$ and $\epsilon>0$ such 
that for $\xi >0$, and $n$ large
\be{UB-prop2}
P\pare{ \sum_{z\in \Z^3} \ind_{\acc{l_n(z)> \log(n)^{\underline\chi}}}
l_n^2(z)> n\xi }\le \exp\pare{-n^{1/3} \log(n)^{\epsilon} }.
\ee
Note that \reff{intro.5} improves \reff{UB-prop2}.
One reason to study $||l_n||_q$ for $q>2$, is that 
the upper bound \reff{intro.4} for $q>2$, yields \reff{intro.5} at once.
More precisely, for $q<q_c(d)$, choose $q'$ with
$q<q'<q_c(d)$, and for any $A>0$, the obvious inequality
\be{intro.9}
\sum_{z\in \Z^d}\ind_{\acc{l_n(z)>A}}
\ l_n^q(z)\le \frac{||l_n||_{q'}^{q'}}{A^{q'-q}},
\ee
implies that
\[
P\pare{ \sum_{z\in \Z^d}\ind_{\acc{l_n(z)>A}}l_n^q(z)\ge 
n\xi} \le
P\pare{ ||l_n||_{q'}^{q'}\ge A^{q'-q} n\xi}.
\]
For $A$ and $n$ large enough, $A^{q'-q}n\xi \ge
2 E[||l_n||_{q'}^{q'}]$. 
Thus, if we set $\beta=\frac{2}{d}\frac{q'-q}{q'-1}>0$, 
then from \reff{intro.4}, we have a constant $c_1(d,q')$ such that
\be{intro.10}
P\pare{ \sum_{z\in \Z^d}\ind_{\acc{l_n(z)>A}}l_n^q(z)\ge  
n\xi} \le \exp\pare{ -c_1(d,q') \xi^{\frac{2}{d}
(\frac{1}{q'-1})}A^\beta n^{1-\frac{2}{d}}}.
\ee
Thus, in order to improve \reff{UB-prop2} in $d=3$,
we were left with studying $q$-fold self-intersections with $2<q<3=q_c(3)$.

In most works on two-fold self-intersection, 
a starting point, which we trace back to the work of 
Westwater~\cite{west} and Le Gall~\cite{legall},
is a decomposition of $||l_n||_2^2$ 
in terms of {\it intersection local times}
of two independent random walks starting at the origin.
However, such a decomposition is restricted to 
$q$-fold self-intersection local times with $q\in \N$:
When $q=2$ and $d=3$
(in the {\it sub-critical regime}) 
Le Gall's decomposition is a first step in obtaining, in \cite{Chen08},
a moderate and large deviations principles.
When $q=3$ and $d\ge 4$ (in the {\it super-critical regime}), 
\cite{FMV} uses a type of Le Gall's decomposition 
to obtain moderate and large deviations estimates.

Here, our starting point is an {\it approximate decomposition} obtained 
by slicing $||l_n||_q^q$ over level sets of the local times, 
for any real $q>1$. This
is based on the following simple inequality. Let $\{b_n,\ n\in \N\}$
be a subdivision of $[1,\infty)$, and
let $l_1$ and $l_2$ be positive
integers (which we think of as the local times of a given
site in each half time-period). Then, for $q>1$, we have the upper bound
\be{intro.8}
(l_1+l_2)^q\le l_1^q+ l_2^q
+2^q\sum_{i=0}^\infty b_{i+1}^{q-2}
\ind_{\acc{b_i \le \max(l_1,l_2) < b_{i+1}}} l_1\times l_2,
\ee
as well as the obvious lower bound: $(l_1+l_2)^q\ge l_1^q+l_2^q$.
The desirable feature of \reff{intro.8} is that on its right hand side,
the $q$-th power of $l_1$ and $l_2$ comes without penalty, whereas
the term $l_1\times l_2$ yields an intersection local times. 
Thus, \reff{intro.8} leads to the following result which plays here the
r\^ole of Le Gall's decomposition of~\cite{legall}.
\bp{prop-legall} For any integers $n$ and $l$, with $2^l<n$, let
$\{n_i,i=1,\dots,2^l\}$ be positive integers summing up to $n$.
Let $\{l_{.}^{(i)},i=1,\dots,2^l\}$ be the local times
of $2^l$ independent random walks starting at 0. 
If $\{b_i,\ i\in \N\}$ is a subdivision of $[1,n]$, then,
\be{low-eq.8}
S_q^{(l)}\le ||l_n||_q^q\le S_q^{(l)}+\sum_{j=1}^{l} \I_j,
\quad\text{where}\quad
S_q^{(l)} \stackrel{\text{law}}{=}\ \sum_{i=1}^{2^l} ||l_{n_i}^{(i)}||_q^q,
\ee
and, for $j=1,\dots,l$, and $m_k=n_{(k-1)2^{l-j}+1}+\dots+n_{k2^{l-j}}$
for $k=1,\dots,2^j$
\be{decomp.5}
\I_j\stackrel{\text{law}}{=}\  
\sum_{k=1}^{2^{j-1}}\sum_{i} 2^qb_{i+1}^{q-1}
\pare{\sum_{z:\ b_i\le  l^{(2k)}_{m_{2k}}(z)< b_{i+1}}
\!\!\! l_{m_{2k-1}}^{(2k-1)}(z)
+\sum_{z:\ b_i\le  l^{(2k-1)}_{m_{2k-1}}(z)< b_{i+1}}
\!\!\! l_{m_{2k}}^{(2k)}(z)}.
\ee
\ep
\br{rem-realq} We first note some natural limitations in using
the approximate decomposition \reff{low-eq.8}. When we deal with 
$\{||l_n||_q^q-E[||l_n||_q^q]\ge \xi_n n\}$ for small $\xi_n$,
we need to bound the difference between $E[||l_n||_q^q]$ and
the expectation of the upper bound in \reff{low-eq.8}. When,
we take $l$ such that $2^l\sim n^{1-\delta_0}$, then this
difference turns out to be of order smaller than
$n^{1-\delta_0/2}$, allowing us to write
\be{why-xibig}
\begin{split}
\acc{||l_n||_q^q-E[||l_n||_q^q]\ge \xi_n n}\subset&
\acc{S_q^{(l)}-E[S_q^{(l)}]\ge \frac{\xi_n}{2}n}\\
&\quad\cup
\acc{\sum_{j=1}^{l} \I_j-E[\I_j]\ge \frac{\xi_n}{2}n-n^{1-\delta_0/2}}.
\end{split}
\ee
\reff{why-xibig} requires that 
$\xi_n\ge n^{-\delta_0/2}$. 
\er
Proposition~\ref{prop-legall} is our initial step in the
proof of Theorems~\ref{theo-d3} and \ref{theo-d4}, and leads to a 
central limit theorem (CLT) for $||l_n||_q^q$ in dimension 4 or more, 
as well as a characterization of the variance of $||l_n||_q^q$.

Before stating results concerning the typical behavior of $||l_n||_q^q$,
we give a heuristic discussion of the proof of Theorem~\ref{theo-d3}
assuming Proposition~\ref{prop-legall}. 
More precisely, we wish to sketch the reasons why
the {\it approximate decomposition} \reff{low-eq.8}, 
reduces large deviations for
the $q$-norm of the local times, to large deviations
for a sum of independent geometric variables. 

Consider a choice of $l$ such that
$2^l$ is close to $n$ in \reff{low-eq.8}, and $n_i\sim n/2^l$. 
Then, $S_q^{(l)}$
is a sum of about $n$ independent terms, each one bounded by
its time-span, $n/2^l$, to the power $q$. 
Recall that the probability of deviating 
from the mean, for a sum of $n$ independent and essentially bounded 
variables, is of order $\exp(-O(n))$ (see Lemma~\ref{lem-LD}
for a precise statement). We can therefore
neglect the contribution of $S_q^{(l)}$ to the excess $q$-norm,
though the mean of $S_q^{(l)}$ is close to $E[||l_n||_q^q]$, as
easily seen from Lemma~\ref{lem-moy} below.
Now, for a fixed
$j$ in $\{1,\dots,l\}$, let $m=n/2^j$, and note that $\I_j$
 is a finite sum of independent terms distributed as
$b^{q-1}\ l_{m}(\tilde \D(b))$, where $\tilde \D(b)$ is 
the set $\{z: \tilde l_{m}(z)\sim b\}$, with $\tilde l_{m}$ 
an independent copy of $l_{m}$, and $b$ spans a subdivision of
$[1,n]$. Since we consider
transient random walks, $l_{m}(\tilde \D(b))$
is bounded by a geometric variable (when fixing $\tilde \D(b)$,
as shown in Lemma 1.2 of \cite{AC04}).
At this point, one normalizes $l_{m}(\tilde \D(b))$.
If $P_0$ and $\tilde P_0$
are the law of the two independent copies of the 
(transient) walk started at 0, define
\be{def-X}
X=\frac{l_{m}(\tilde \D(b))}{E_0[l_{m}(\tilde \D(b))]},
\quad\text{ so that for some }\kappa>0,\qquad
\tilde P_0\otimes P_0(X>t)\le e^{-\kappa t}.
\ee
Now, it is well known that for any $m$,
$E_0[l_{m}(\tilde \D(b))]\le C |\tilde \D(b)|^{2/d}$, (see
Lemma~\ref{lem-app.1}). Thus, an estimate of the large deviation
probability requires an estimate on the volume of level sets of the
local times. Now, in obtaining
a bound on the volume of $\tilde \D(b)$, assume for simplicity
that we only have two types of $b$: that is,
we distinguish {\it often visited sites}, say sites visited $n^x$-times
with $x$ close to $1/q$,
whose level sets are part of what we call {\it top levels}, and say the 
once-visited sites, whose level sets are part of what we
call {\it bottom levels}.
The {\it bottom levels} are the easiest to 
treat (see Section~\ref{sec-bottom}, and Lemma~\ref{lem-bottom}). 
Indeed, we use essentially that $|\tilde \D(b)|\le n$ for $b\sim 1$,
so that we expect (when we restrict $\I_j$ only to {\it bottom levels}
and using $X_k=X$ in law)
\[
P\pare{\I_j-E[\I_j]\ge n\xi}\sim P\pare{
\sum_{k=1}^{2^j} X_k-E[X_k]\ge \frac{n\xi}{n^{2/d}}}\sim 
\exp\pare{-O\pare{n^{1-\frac{2}{d}}}}.
\]

The {\it top levels}, treated in Section~\ref{sec-top},
require a type of {\it bootstrap argument}, using that if 
$\tilde \D(b)$ has a large volume, it 
implies that the $q$-norm of the local time 
$\tilde l_{m}$ is large. The bootstrap argument yields
Corollary~\ref{cor-legall}. It allows us to normalize properly
our geometric-like random variables $b^{q-1}\ l_m(\tilde \D(b))$.

We turn now to the typical behavior of the q-norm. Chen has provided
in~\cite{Chen07} asymptotics for the variance 
of $||l_n||^2_2$ in $d\ge 3$. He shows that
(i) in $d=3$, $\var(||l_n||^2_2)\sim \lambda_3 n\log(n)$, and
(ii) in $d\ge 4$, $\var(||l_n||^2_2)\sim \lambda_d n$, where $\lambda_d$ are
constants expressed in terms of the Green's function of the walk. Following
ideas of Jain and Pruitt~\cite{jainpruitt}, and of
Le Gall and Rosen~\cite{legallrosen},
Chen obtains a CLT in dimension 3 or more for $||l_n||^2_2$. 
Finally, Becker and K\"onig in \cite{BK} have shown that for $q$ integer, 
(i) in $d=3$, $\var(||l_n||_q^q)\le n^{3/2}$, (ii) in $d=4$, $\var(||l_n||_q^q)\le 
n\log(n)$, and (iii) in $d\ge 5$, $\var(||l_n||_q^q)\le c_d n$. 
Our result deals with the general case ($q>1$ real), 
where no representation of $||l_n||_q^q$ is possible
in terms of multiple time-intersections. We transform
Lindeberg's condition into a {\it large deviation} event for 
$||l_{T_n}||^2_2$ on
the scale of time of the CLT, that is $T_n\approx\sqrt{n}$.

We start with an estimate for the expectation of $||l_n||_q^q$,
of the same type as
Theorem 1 of Dvoretzky and Erd\"os~\cite{DE} for the range of a transient random walk.
Thus, if $\gamma_d$ is the 
probability of never returning to its original position, it is shown
in~\cite{DE} that for positive constants $c_d$, when $R_n$ 
is the set of visited sites before time $n$,
\be{moy.1}
|E[|R_n|]-n\gamma_d|\le c_{d} \psi_d(n), \quad\text{with}\quad
\psi_d(n)=
\left\{ \begin{array}{ll}
n^{1/2} & \mbox{ for } d=3 \, , \\
\log(n) & \mbox{ for } d=4 \, , \\
1 & \mbox{ for } d\ge 5 \, ,
\end{array} \right.
\ee
Jain and Pruitt~\cite{jainpruitt}
obtain the asymptotics $\var(|R_n|)\sim a\log(n)n$ for some $a>0$ in 
$d=3$, and $\var(|R_n|)\sim c_d' n$ in $d>3$,
for some positive constants $c_d'$. The corresponding CLT (in
$d\ge 3$) was shown by
Jain and Pruitt~\cite{jainpruitt} for the simple random walk, 
and by Le Gall and Rosen~\cite{legallrosen} for stable random walks.
Note that the limiting law is gaussian, in $d\ge 3$ but fails to 
be so in $d=2$, as shown by Le Gall in~\cite{legall86}.
\bl{lem-moy}
Assume that $d\ge 3$ and $q>1$. There is a constant $C_d$, such that
\be{moy-key}
0\le \kappa(q,d) n- E[||l_n||_q^q]\le C_d\psi_d(n), \quad\text{with}\quad
\kappa(q,d)=\gamma_d E[l_{\infty}(0)^q].
\ee
Also, if $d=3$, then, there is a constant $c_3$ such that
\be{key-var1}
\var(||l_n||_q^q)\le c_3 \log(n)^2\ n.
\ee
If $d\ge 4$, then there are positive constants $v(q,d)$ and $c(q,d)$,
such that
\be{key-var2}
|\frac{\var(||l_n||_q^q)}{n}-v(q,d)|\le c(q,d) \frac{\log(n)}{\sqrt{n}}.
\ee
\el
Finally, we have the following central limit theorem.
\bt{theo-clt}
If $Z$ is a standard normal variable, then
\be{key-clt}
\frac{||l_n||_q^q-n\kappa(q,d)}{\sqrt{n v(q,d)}}
\quad\stackrel{\text{law}}{\longrightarrow}\quad Z.
\ee
\et
A challenging open question is to to understand the strategy which
realizes $\{||l_n||_q^q-E[||l_n||_q^q]\ge \xi n\}$, right at
the critical value $q=q_c(d)=\frac{d}{d-2}$.

The paper is organized as follows. 
The approximate decomposition of
$||l_n||_q^q$ is given in Section~\ref{sec-decomp}.
The sub-critical regime is studied in Section~\ref{sec-d3}:
The upper bound in \reff{intro.4} is proved in Section~\ref{sec-proofd3},
and the lower bound is given in Section~\ref{sec-LBd3}. The super-critical
regime is studied in Section~\ref{sec-supercritical}. 
Theorem~\ref{theo-d4} is proved in Section~\ref{sec-supercritical}.
The proof of \reff{intro.7} is given in Section~\ref{sec-superlow}.
The proof of Lemma~\ref{lem-level} is given in Section~\ref{sec-gamma}. 
Lemma~\ref{lem-moy}, as well as the CLT are proved in Section~\ref{sec-clt}.
In Section~\ref{sec-append}, we recall Lemma 5.1 of \cite{A06}, and
improve Lemma 5.3 of \cite{A06}, used to
control intersection local times-type quantities.

\section{General Considerations ($q>1$)}\label{sec-general}
In this section, we deal with the general case $q>1$.
In section \ref{sec-decomp}, we develop a approximation of $||l_n||_q^q$ as
sums of two types of independent variables.
\begin{enumerate}
\item Intersection local times of independent walks.
\item Self-intersection local times, on a much shorter time-period.
\end{enumerate}
In section \ref{sec-self}, we treat the sums of
self-intersection local times. 
\subsection{Approximate decomposition for $||l_n||_q^q$}\label{sec-decomp}
Before we prove 
Proposition~\ref{prop-legall}, we present a useful corollary
which requires more notations.

For integers $n$ and $l$, with $2^l<n$, we 
recall the ``almost'' dyadic decomposition of $n$ of
Remark 2.1 of~\cite{AC05}. We divide $n$ into $2^l$ integers
$n_1^{(l)},\dots,
n_{2^l}^{(l)}$ with $n=n_1^{(l)}+\dots+n_{2^l}^{(l)}$ and
\be{quasi-dyadic}
\max_i(n_i^{(l)})-\min_i(n_i^{(l)})\le 1,\quad
\frac{n}{2^l}-1\le n_i^{(l)}\le \frac{n}{2^l}+1,
\quad\text{and}\quad
n_{k}^{(l-1)}=n_{2k-1}^{(l)}+n_{2k}^{(l)}.
\ee
We run $2^l$ independent random walks starting at the origin. 
The $i$-th walk
runs for a time-period $[0,n_i^{(l)}[$, and we denote by 
$l_i^{(l)}:\Z^d\to\N$
its local times during time-period $[0,n_i^{(l)}[$. 
Also, we introduce, for $k=1,\dots,2^{l}$, the following sets
\be{turb-5}
\D^{(l)}_{k,i}=\acc{z\in \Z^d:\ b_i\le  l^{(l)}_k(z)< b_{i+1}}.
\ee
Now, for any $M>0$, let $\{b_i,\ i\in \N\}$ be a subdivision 
of $[1,M]$, and denote by $\Theta_M(x)=x\ind_{\{x\le M\}}$.
\br{rem-legall2}
We could restrict the sum over $\Z^d$ which enters $||l_n||_q^q$ in
\reff{low-eq.8} over $\{z:\ l_n(z)\le M\}$ for any positive $M$.
The proof of Proposition~\ref{prop-legall} yields, for
any $\{b_i,\ i\in \N\}$ subdivision of $[1,M]$, 
\be{low-levelq}
\sum_{z\in \Z^d} \Theta_M\pare{l_n(z)}^q\le  
\sum_{k=1}^{2^l}\sum_{z\in \Z^d} \Theta_M\pare{l_k^{(l)}(z)}^q+
\sum_{j=1}^{l} \I_j.
\ee
The only difference with \reff{low-eq.8} is the subdivision
which enters into the definition of $\I_j$.
The proof of Proposition~\ref{prop-legall} 
is written in view of \reff{low-levelq} (see the key step \reff{decomp.2}).
\er
As a corollary of \ref{low-levelq}, we obtain the following result.
\bc{cor-legall} For any $M>0$,
let $\{b_i,\ i\in \N\}$ be a subdivision of $[1,M]$.
For any integers $n$ and $L$, with $2^L<n$, and
for any sequence of positive numbers $\{m_n,\epsilon_n,\ n\in \N\}$, 
we have
\be{bootstrap.2}
\begin{split}
P\pare{ ||\Theta_M(l_n)||_q^q\ge m_n+\epsilon_n}\le &
2^{L+1} P\pare{ \sum_{j=1}^{2^L} ||\Theta_M(l_j^{(L)})||_q^q\ge m_n}\\
&+\sum_{h=1}^{L} 2^{h} P\pare{ \sum_{l=h}^{L}
\ind_{\acc{\G^{(l)}_1\cap \G^{(l)}_2}}\I_l\ge \epsilon_n},
\end{split}
\ee
where for $l\le L$, $k=1,\dots,2^{l}$, and $i\in \N$
\be{bulk.1}
\G^{(l)}_{k,i}=\acc{|\D^{(l)}_{k,i}|\le \frac{m_n+\epsilon_n}{b_i^q }},
\quad \G^{(l)}_1=\bigcap_{k=1}^{2^{l}}\bigcap_i \G^{(l)}_{2k,i}
\quad\text{and}\quad 
\G^{(l)}_2=\bigcap_{k=1}^{2^{l}}\bigcap_i \G^{(l)}_{2k-1,i}.
\ee
\ec
\br{rem-legall1}
The symbols, $\epsilon_n$ and $m_n$, are suggestive of the fact that
when $L$ is large enough, the sum of $2^L$
independent $q$-fold self-intersections, that we called $S_q^{(L)}$,
stays close to its mean, which is also close to the mean of $||l_n||_q^q$.
This is shown in Section~\ref{sec-self}. So, $m_n$ stands for {\it mean},
and $\epsilon_n$ stands for {\it excess}. To estimate how small can
$\epsilon_n$ be, we now compute the expectation of $\sum_{l=1}^L\I_l$.
We use Lemma~\ref{lem-app.1} in the worse case, that is in dimension 3,
to obtain for some constants $c_1,c_2$ and $c_3$
\ba{bulk-center2}
E\cro{ \sum_{l=1}^L\I_l} &=& 
\sum_{l=1}^{L}2^l\sum_{i\in \N}
2^{q}(b_{i+1})^{q-1}C_d\psi_d(\frac{n}{2^l})
e^{-\kappa_d  b_i}\cr
&\le & c_1{\sqrt n}  \sum_{l=1}^{L} 2^{l/2} \sum_{i\in \N}
(b_{i+1})^{q-1} e^{-\kappa_d b_i}\cr
&\le & c_2 \sqrt{ 2^{L}n} \sum_{i\in \N}
(b_{i+1})^{q-1} e^{-\kappa_d b_i} .
\ea
First, we need to choose a subdivision $\{b_i,i\in\N\}$ such that
the last sum in \reff{bulk-center2} is convergent. 
Secondly, the right hand side of \reff{bootstrap.2} is small if
$ \sum_{l=h}^{L}E[ \ind\{\G^{(l)}_1\cap \G^{(l)}_2\}
\I_l]\ll \epsilon_n$. From \reff{bulk-center2}, we see that
$\epsilon_n$ can be chosen small if $2^L\ll n$. 
On the other hand, we see in Section~\ref{sec-self},
that $L$ has to be large enough,
for the probability of $\{S_q^{(L)}\ge m_n\}$ to be negligible, 
when $m_n=E[||l_n||_q^q]+n\epsilon$. Remark~\ref{rem-self1} shows
that a choice of $L$ such that
$2^L>n^{1-\delta_0}$ with $q \delta_0<\frac{2}{d}$, fulfills both 
requirements.

\er

\noindent{\bf Proof of Proposition~\ref{prop-legall}}

The proof proceeds by induction on $l\ge 1$. It is however
easy to see that proving the case $l=1$ requires the same
arguments as going from $l-1$ to $l$. We focus on 
the first step $l=1$, and omit the easy passage from $l-1$ to $l$.

For any $x\in [0,1]$, and $q\ge 1$, we have
\be{low-eq.0}
1+x^q\le (1+x)^q\le 1 + x^q+2^q x.
\ee
Thus, for any nonnegative
integers $l_1,l_2$ with $0\le l_1,l_2\le M$, we have from \reff{low-eq.0}
\be{low-eq.1}
l_1^q+l_2^q\le (l_1+l_2)^q\le l_1^q+l_2^q+2^q M^{q-2} l_1l_2.
\ee
Now, for any $M>0$, let $\{b_i,\ i\in \N\}$ be a subdivision 
of $[1,M]$, and recall that $\Theta_M(x)=x\ind_{\{x\le M\}}$.
For any nonnegative integers $l_1,l_2$ 
\be{low-eq.2}
\begin{split}
\pare{ \Theta_M(l_1+l_2)}^q\le&\ \pare{\Theta_M(l_1)}^q
+ \pare{ \Theta_M(l_2)}^q\\
&\quad+2^q\sum_{i=1}^n b_{i+1}^{q-2}
\ind\acc{b_i\le \max(l_1,l_2) < b_{i+1}}\ l_1\times l_2.
\end{split}
\ee 
Indeed, $l_1+l_2\le M$ and $l_1,l_2\ge 0$, imply (i) $l_1\le M$ and $l_2
\le M$,
and (ii) for some $i_0>0$, $\max(l_1,l_2)\in[b_{i_0},b_{i_0+1}[$. 
Then, from \reff{low-eq.1}
\[
\pare{ \Theta_M(l_1+l_2)}^q\le 
\Theta_M\pare{ l_1}^q+
\Theta_M\pare{ l_2}^q+ 2^qb_{i_0+1}^{q-2}l_1\times l_2. 
\]
For any integer $n$, we consider the local time $l_n$, which
we denote as $l_{[0,n[}(z)$ to emphasize the time period. For any
integer $n_1$ with $0<n_1<n$, set $n_2=n-n_1$, and 
from the increments of our initial random walk, say $\{Y_n,\ n\in \N\}$,
we build two independent random walks with local times 
\[
l^{(1,1)}_{]0,k]}(z)=\ind\{Y_{n_1}=z\}+\dots+
\ind\{Y_{n_1}+\dots+Y_{n_1-k+1}=z\},
\]
and,
\be{eq-not.1}
l^{(1,2)}_{[0,k[}(z)=\ind\{0=z\}+\ind\{-Y_{n_1+1}=z\}+\dots+
\ind\{-Y_{n_1}-\dots-Y_{n_1+k-1}=z\}.
\ee
It is obvious that on $\{S(n_1)=y\}$, 
\be{decomp.1}
l_n(z)=l^{(1,1)}_{]0,n_1]}(y-z)+l^{(1,2)}_{[0,n_2[}(y-z).
\ee
If we denote
$\bar l^{(1)}(z)=\max(l^{(1,1)}_{]0,n_1]}(z), l^{(1,2)}_{[0,n_2[}(z))$,
and sum \reff{decomp.1} over $z\in \Z^d$, we obtain 
\be{decomp.2}
\begin{split}
\sum_{z} &\Theta_M\pare{l_n(z)}^q\le \sum_{z}
\Theta_M\pare{l^{(1,1)}_{]0,n_1]}(S(n_1)-z)}^q+  \sum_{z}
\Theta_M\pare{l^{(1,2)}_{[0,n_2[}(S(n_1)-z)}^q\\
&+2^q\sum_{z\in \Z^d} \sum_{i=1}^n b_{i+1}^{q-2}
\ind_{\acc{b_i\le \bar l^{(1)}(S(n_1)-z) < b_{i+1}}}\
l^{(1,1)}_{]0,n_1]}(S(n_1)-z)\times l^{(1,2)}_{[0,n_2[}(S(n_1)-z)\\
&\le \sum_{z\in \Z^d} \Theta_M\pare{l^{(1,1)}_{]0,n_1]}(z)}^q+
\ \sum_{z\in \Z^d} \Theta_M\pare{l^{(1,2)}_{[0,n_2[}(z)}^q\\
&+ 2^q\sum_{z\in \Z^d}\sum_{i=1}^n b_{i+1}^{q-2}
\ind_{\acc{b_i\le \bar l^{(1)}(z) < b_{i+1}}}\
l^{(1,1)}_{]0,n_1]}(z)\times l^{(1,2)}_{[0,n_2[}(z).
\end{split}
\ee
Now, we rewrite \reff{decomp.2} in a concise form as
\be{decomp.3}
||\Theta_M( l_n )||_q^q\le 
||\Theta_M( l^{(1,1)}_{]0,n_1]})||_q^q+
||\Theta_M( l^{(1,2)}_{[0,n_2[)})||_q^q+\tilde\I_1(n_1,n_2),
\ee
where the term dealing with intersection times of independent strands is
\be{decomp.4}
\tilde\I_1(n_1,n_2)=2^q\sum_{z\in \Z^d}\sum_{i=1}^n b_{i+1}^{q-2}
\ind_{\acc{b_i\le \bar l^{(1)}(z) < b_{i+1}}}\
l^{(1,1)}_{]0,n_1]}(z)\times l^{(1,2)}_{[0,n_2[}(z).
\ee
To prove the upper bound in \reff{low-eq.8} for $l=1$, it suffices to set
$M=n$ (so that this truncation disappears), and to show that
$\tilde\I_1(n_1,n_2)\le \I_1$ given in \reff{decomp.5}. This latter
inequality follows first by noting the obvious inclusion
\[
\begin{split}
\acc{z:\ b_i\le \max(l^{(1,1)}_{]0,n_1]}(z),l^{(1,2)}_{[0,n_2[}(z))
< b_{i+1}} \subset& \acc{z:\ b_i\le  l^{(1,1)}_{]0,n_1]}(z)< b_{i+1}}\\
&\cup \acc{z:\ b_i\le  l^{(1,2)}_{[0,n_2[}(z)< b_{i+1}},
\end{split}
\]
and secondly, using that
\be{decomp.7}
\begin{split}
\sum_{z\in \Z^d} 
\ind_{\acc{b_i\le \bar l^{(1)}(z) < b_{i+1}}}\
l^{(1,1)}_{]0,n_1]}(z)\times l^{(1,2)}_{[0,n_2[}(z)
&\le
\sum_{z\in \Z^d}\ind\acc{b_i\le l^{(1,1)}_{]0,n_1]}(z)< b_{i+1}}b_{i+1} 
l^{(1,2)}_{[0,n_2[}(z)\\
&+ \sum_{z\in \Z^d}\ind\acc{b_i\le  l^{(1,2)}_{[0,n_2[}(z)< b_{i+1}}b_{i+1} 
l^{(1,1)}_{]0,n_1]}(z).
\end{split}
\ee
As we iterate the approximate decomposition $l$-times, we 
obtain the upper bound in \reff{low-eq.8}, or more generally
the bound \reff{low-levelq}.

The lower bound in \reff{low-eq.8} is an obvious corollary of the
inequality $(l_1+l_2)^q\ge l_1^q+l_2^q$, valid for $q\ge 1$ and
$l_1,l_2$ nonnegative.

\qed

\noindent{\bf Proof of Corollary~\ref{cor-legall}}

We write the case $M=n$, that is the case with no truncation. The
case with truncation is obtained as we replace $l_n(z)$ by $
\Theta_M(l_n(z))$ wherever it appears.
Assume that we stop the induction in Proposition~\ref{prop-legall} 
at some step $L$ (typically $2^L=n^{1-\delta_0}$ and $\delta_0$ small).
For any sequence of positive numbers $\epsilon_n,m_n$, 
we have from \reff{low-eq.8},
\be{low-key1}
P\pare{  ||l_n||_q^q \ge m_n+\epsilon_n}\le
P\pare{ S_q^{(L)}\ge m_n}+ P\pare{\sum_{l=1}^{L} \I_l\ge \epsilon_n},
\text{  where  }
S_q^{(l)}=\sum_{k=1}^{2^l} ||l_{k}^{(l)}||_q^q. 
\ee
We introduce, as in \cite{AC05,A06}, a bootstrap control on the
volume of $D^{(l)}_{k,i}$. Consider $\G^{(l)}_{k,i}$
given in \reff{bulk.1}. 
On the complement of $\G^{(l)}=\G^{(l)}_1\cap \G^{(l)}_2$, 
there is $k_0,i_0$ such that
$|D^{(l)}_{k_0,i_0}|>(m_n+\epsilon_n)/b_{i_0}^q$ so that
\be{bulk.2}
S_q^{(l)}\ge \sum_{z\in D^{(l)}_{k_0,i_0}}\pare{ l^{(l)}_{k_0}(z)}^q
\ge \frac{m_n+\epsilon_n}{b_{i_0}^q} b_{i_0}^q=m_n+\epsilon_n.
\ee
Writing $S_q^{(0)}=||l_n||_q^q$, we write a more suggestive relation
\be{bootstrap.1}
\begin{split}
P\pare{ S_q^{(0)}\ge m_n+\epsilon_n}\le & 
P\pare{ S_q^{(L)}\ge m_n}+
P\pare{\sum_{l=1}^{L} \ind_{\acc{\G^{(l)}}}\I_l\ge \epsilon_n}\\
&+\qquad \sum_{l=1}^{L}P\pare{ S_q^{(l)}\ge m_n+\epsilon_n}.
\end{split}
\ee
Starting the approximation with $S_q^{(l)}$ with $l<L$, we obtain similarly
\be{induc.1}
P\pare{ S_q^{(l)}\ge m_n+\epsilon_n}\le P\pare{ S_q^{(L)}\ge m_n}+
P\pare{\sum_{j=l+1}^{L}\!\! \ind_{\acc{\G^{(j)}}}\I_j\ge \epsilon_n}
+\sum_{j=l+1}^{L}\!\!P\pare{ S_q^{(j)}\ge m_n+\epsilon_n}.
\ee
Assume now that for $j>l$, and $j<L$ we have
\be{induc.2}
P\pare{ S_q^{(j)}\ge m_n+\epsilon_n}\le 2^{L-j+1} P\pare{S_q^{(L)}\ge m_n}
+\sum_{h=j+1}^{L} 2^{h-j-1} P\pare{ \sum_{i=h}^{L}\ind_{\acc{\G^{(i)}}}
\I_i\ge \epsilon_n}.
\ee
Note that \reff{induc.2} is true for $j=L-1$. 
Then, using the hypothesis \reff{induc.1} in \reff{induc.2}, we obtain
\be{induc.3}
\begin{split}
P\pare{ S_q^{(l)}\ge m_n+\epsilon_n}\le &
\sum_{j=l}^L 2^{L-j} P\pare{ S_q^{(L)}\ge m_n}
+\sum_{j=l+1}^{L} \sum_{h=j+1}^{L} 2^{h-j-1} 
P\pare{\sum_{i=h}^{L}\ind_{\acc{\G^{(i)}}}\I_i\ge \epsilon_n}\\
\le & 2^{L-l+1} P\pare{ S_q^{(L)}\ge m_n}+
\sum_{h=l+1}^{L}2^h\pare{\sum_{l<j<h}2^{-j}} 
P\pare{\sum_{i=h}^{L}\ind_{\acc{\G^{(i)}}}\I_i\ge \epsilon_n}
\end{split}
\ee
By way of induction, \reff{induc.3} yields \reff{bootstrap.2}.
\qed
\subsection{On large sums of $q$-fold self-intersection}
\label{sec-self}
In this section, we consider the contribution of the term $S_q^{(l)}$,
which appears in \reff{low-eq.8}, in making $\{||l_n||_q^q-E[||l_n||_q^q]\ge
\xi n\}$. Recall that $S_q^{(l)}$ is a sum of independent
copies of $q$-fold self-intersection over times of order $n/2^l$.
We first choose $l$ large enough, and
then use the boundedness of the $q$-fold self-intersection.

Fix $\delta_0$ such that $0<\delta_0<\frac{2}{qd}$. Let $L$ be an
integer so that $2^{L}\le n^{1-\delta_0}<2^{L+1}$. Note the obvious bound
\be{supLq}
\max_{k\le 2^L} ||l_{k}^{(L)}||_q^q\le
\max_{k\le 2^L} (n_k^{(L)})^{q}\le \pare{ \frac{n}{2^L}+1}^q\le
2^{q+1} n^{q\delta_0}.
\ee
The main result, in this section, reads as follows.
\bl{lem-LD}
Fix $\delta\ge 0$, with either (i) dimension is 3 and $\delta< \delta_0/2$,
or (ii) dimension is 4 or more and $\delta<\delta_0$.
Let $\xi_n\ge n^{-\delta}$. Then, for $n$ large enough
\be{self-key1}
P\pare{ |S_q^{(L)}-E[||l_n||_q^q]|> 
\xi_n n}\le \exp\pare{-\frac{\xi_n}{2} n^{1-\delta_0 q}}
\ee
\el
\br{rem-self1}
Let us consider now the regimes of Theorems~\ref{theo-d3} and
\ref{theo-d4}.
\begin{itemize}
\item When $q<q_c(d)$, the speed exponent in \reff{intro.4} is
$1-\frac{2}{d}$. Thus, the right hand side of \reff{self-key1} 
with $\xi_n=\xi$ is negligible when $1-q\delta_0>1-\frac{2}{d}$, 
so that we need $q\delta_0<2/d$.
\item When $q>q_c(d)$, the speed exponent in \reff{intro.4} is
$\frac{1}{q}$. It is enough to have again $q\delta_0<2/d$.
\end{itemize}
\er

\noindent{{\bf Proof of Lemma~\ref{lem-LD}}}

First, we write
\be{sup-eq1}
S_q^{(L)}-E[||l_n||_q^q]=\sum_{k=1}^{2^L} Z(k)\ +R_1,\quad\text{with}\quad
Z(k)=||l_{k}^{(L)}||_q^q-E\cro{||l_{k}^{(L)}||_q^q},
\ee
and
\[
R_1=\sum_{k=1}^{2^L}\pare{E\cro{||l_{k}^{(L)}||_q^q}
-\kappa(q,d)n_k^{(L)}}\quad-\quad\pare{E[||l_n||_q^q]-\kappa(q,d)n}.
\]
Using Lemma~\ref{lem-moy} in $d\ge 3$, we have a constants $c_d$ such that
\be{sup-eq2}
|R_1|\le c_d\psi_d(n)+c_d\sum_{k=1}^{2^L}\psi_d\pare{n_k^{(L)}}\le 
c_d\pare{\psi_d(n)+ 2^{L+1}\psi_d\pare{\frac{n}{2^L}}}.
\ee
Thus, for $\xi_n\ge n^{-\delta}$ and (i) $0\le \delta<\delta_0/2$
in $d=3$, or (ii) $0\le \delta<\delta_0$ in $d>3$, we have
\be{sup-eq3}
P\pare{|S_q^{(L)}-E[||l_n||_q^q]|\ge \xi_n n}\le 
P\pare{|\sum_{k=1}^{2^L}Z(k)|\ge \frac{\xi_n}{2} n}.
\ee
We note that $|x|=\max(x,-x)$, and use Chebyshev's exponential inequality.
For $\lambda\in [0,1]$, 
\be{sup-eq4}
\begin{split}
P\pare{\pm\sum_{k=1}^{2^L} Z(k)\ge \frac{\xi_n}{2} n}&\le
\exp\pare{ -\lambda \frac{\xi_n}{2} (\frac{2^L}{n})^q n}
\pare{E\cro{ e^{\pm\lambda (\frac{2^L}{n})^q Z(k)}}}^{2^L}\\
&\le \exp\pare{ -\lambda \frac{\xi_n}{2} (\frac{2^L}{n})^q n}
\pare{1+\lambda^2 (\frac{2^L}{n})^{2q}\var(Z(k))}^{2^L}\\
&\le \exp\pare{ -\lambda \frac{\xi_n}{2} (\frac{2^L}{n})^q n+
\lambda^2 2^L (\frac{2^L}{n})^{2q}\var(Z(k))}.
\end{split}
\ee
We used the uniform bound \reff{supLq} on $|Z(k)|$ in the second 
inequality, and the fact that for $x\le 1$, we have $e^x\le 1+x+x^2$. 
We recall that the bound \reff{key-var1} holds in dimension 3
and more, and reads $\var(Z(k))\le \frac{n}{2^L}
\log^2(\frac{n}{2^L})$. Thus, \reff{sup-eq4} is useful if
\be{sup-eq5}
\frac{\xi_n}{2} (\frac{2^L}{n})^q n
\ge 2\lambda 2^L \frac{n}{2^L}\log^2(\frac{n}{2^L})
(\frac{2^L}{n})^{2q}
\Longleftrightarrow
\xi_n\ge 4\lambda \log^2(\frac{n}{2^L}) (\frac{2^L}{n})^{q}.
\ee
Since $\xi_n\ge n^{-\delta}$, 
\reff{sup-eq5} is implied if $\delta_0 q>\delta$,
which holds if conditions (i) or (ii) of Lemma~\ref{lem-LD} are assumed.

In case (i) or (ii), we choose $\lambda=1$, and take
$n$ large enough so that \reff{sup-eq5} holds. We then obtain
\reff{self-key1}.
\section{The sub-critical regime}\label{sec-d3}
We consider here the case $q<\frac{d}{d-2}$.
The main result of this section is the upper bound of \reff{intro.4}. 
Indeed, we have shown in the Introduction (in \reff{intro.10}) that
\reff{intro.4} implies \reff{intro.5}. 
Finally, the easy lower bound in \reff{intro.4} 
is proved in Section~\ref{sec-LBd3}. 

We have divided the proof into many sections. 
Our starting point is \reff{low-eq.8}.
With the notation of Section~\ref{sec-decomp}, we set,
with $2\epsilon_0<1$,
\[
m_n=E[||l_n||_q^q]+n\epsilon_0\xi_n,\quad\text{ and }
\quad \epsilon_n=n\xi_n(1-\epsilon_0).
\]
In Section~\ref{sec-subdivision}, we choose an appropriate subdivision
of $[1,n]$. When $q<q_c(d)$, strategy B described in the Introduction
suggests to divide the set of visited sites into those
visited about $\xi_n^{1/(q-1)}$-times, and the remaining 
{\it too often visited sites}. 
The contribution of the former sites to $\sum \I_l$ in \reff{low-eq.8},
is called the bottom-level term, and is treated in
Section~\ref{sec-bottom}. The contribution of the 
latter sites is called the
top-level term, and is treated in Section~\ref{sec-top}. 

\subsection{A choice of a subdivision}\label{sec-subdivision}
We first choose the largest $\alpha_0$ such that
\be{alpha0}
\alpha_0^{q-1} \sum_{l=0}^{\infty} 
\pare{ \frac{1}{2^{q-1}}}^{l}\le \frac{1}{16},
\ee
and, for some positive integer $j_0$, $\alpha_0 \xi_n^\gamma=2^{j_0}$. 
Note that $\alpha_0$ is bounded by $1$, though $j_0$ grows with $\xi_n$. 
Recall that $\gamma=\frac{1}{q-1}$, and consider for $i=-j_0,\dots,M_n$
\be{bitop}
b_i=\xi_n^{\gamma} \beta_i,
\quad\text{with}\quad \beta_i=\alpha_0 2^i,
\ee
where $M_n$ is such that $\beta_{M_n}$ is of order $n^{1/q_c(d)}$.
We divide the intersection local times according to 
whether $l^{(l)}_{k}(z)\ge \alpha_0\xi_n^\gamma$ 
(which yields what we call a top-level term),
or $l^{(l)}_{k}(z)< \alpha_0\xi_n^\gamma$ (which yields what 
we call a bottom-level term). Introduce, for $l\le L$
\be{low-eq.10}
\C_n^{\uparrow}(l)=\sum_{i\ge 0}\sum_{k=1}^{2^{l}}
2^q(\xi_n^\gamma \beta_{i+1})^{q-1}\pare{\ind_{\G^{(l)}_{2k,i}}
l_{2k-1}^{(l)} \pare{\D^{(l)}_{2k,i}}+
\ind_{\G^{(l)}_{2k-1,i}} l_{2k}^{(l)}\pare{\D^{(l)}_{2k-1,i}}},
\ee
where for a subset $A$ of $\Z^d$, $l_k^{(l)}(A)=\sum_{z\in A}
l^{(l)}_k(z)$ and,
\be{low-eq.11}
\C_n^{\downarrow}(l)=\sum_{-j_0\le j< 0}\sum_{k=1}^{2^{l}}
2^q(\xi_n^\gamma \beta_{j+1})^{q-1}\pare{
l_{2k-1}^{(l)}\pare{\D^{(l)}_{2k,j}}+
l_{2k}^{(l)}\pare{\D^{(l)}_{2k-1,j}}}.
\ee
Note that if for any $\alpha_0$ satisfying \reff{alpha0}, 
we have $\alpha_0 \xi_n^\gamma<1$,
then there will be no term $\C^\downarrow_n(l)$. 
Note also that in both $\C_n^{\uparrow}(l)$ and $\C_n^{\downarrow}(l)$, the 
sum over $k$ is over independent variables. 
We call $\C_n^{\uparrow}(l)$ the top-level
term, and $\C_n^{\downarrow}(l)$ the bottom-level term.

We choose now $L$ such that $2^L=n^{1-\delta_0}$, and
inequality \reff{bulk-center2} of Remark~\ref{rem-legall1} 
gives us that for some constant $c_3$
\be{low-center}
\sum_{\dagger\in \{\uparrow,\downarrow\}} 
\sum_{h=1}^L E[\C^\dagger(h)]\le c_3 n^{1-\frac{\delta_0}{2}}.
\ee
We denote by $\bar \C^\dagger(h)=\C^\dagger(h)-E[\C^\dagger(h)]$.
Finally, $\xi_n n\ge 4 c_3n^{1-\delta_0/2}$ implies that 
$\xi_n\ge 4c_3n^{-\delta_0/2}$. Inequality \reff{bootstrap.2} yields 
\be{bulk.3}
\begin{split}
P\pare{ ||l_n||_q^q-E[||l_n||_q^q]\ge n\xi_n}\le& 
2^L P\pare{ S_q^{(L)}-E[||l_n||_q^q]\ge n \epsilon_0}\\
&+\sum_{h=1}^{L-1} 2^{h-1}\sum_{\dagger\in \{\uparrow,\downarrow\}}
P\pare{ \bar \C^\dagger(h)\ge \frac{n\xi_n}{8}}.
\end{split}
\ee
Note that from Lemma~\ref{lem-LD}, we have
\be{cor-LD}
P\pare{ S_q^{(L)}-E[||l_n||_q^q]\ge \epsilon_0 \xi_n n}
\le \exp\pare{-\frac{\epsilon_0\xi_n}{2} n^{1-\delta_0 q}}.
\ee
\reff{cor-LD} shows that the contribution of $S_q^{(L)}$ to an excess
self-intersection local times is negligible when $1-\delta_0 q>1-
\frac{2}{d}$, that is when $q\delta_0<\frac{2}{d}$. It remains to show
that the other terms in \reff{bulk.3} are of the right order. 
\subsection{The {\it bottom-level} terms }\label{sec-bottom}
Note that the bottom-level sets $\C_n^\downarrow$ depend on
$\xi_n$. Also, from \reff{cor-LD}, we need only consider 
generation $l<L$ with $2^L=n^{1-\delta_0}$ for $q\delta_0<\frac{2}{d}$.
We establish in this section, the following result.
\bl{lem-bottom} Assume $d\ge 3$ and $q>1$. There is a constant $C>0$
such that for any $h\in \{1,\dots,L\}$, and $1\le \xi_n$
\be{bottom.1}
P\pare{ \sum_{l=h}^L \bar \C_n^\downarrow(l)\ge \frac{\xi_n n}{8}}\le
\exp\pare{ -C \xi_n^{\frac{2}{d}\frac{1}{q-1}}n^{1-\frac{2}{d}}}.
\ee
\el
\br{rem-bottom} 
Recall that $\alpha_0\le 1$, and that if $\xi_n<1$, then 
the terms $\{\C_n^\downarrow(l)\}$ vanish.
\er
\bpr
We first show that we can restrict the sum over $i$ in the definition of
$\C_n^\downarrow(l)$ in \reff{low-eq.11}, to $0>i\ge-l$. 
We make use of the obvious fact
that for any generation $l$, the total time over which run the local times
of the $2^l$ strands is $n$. In other words,
\[
\sum_{k=1}^{2^l} \sum_{z\in \Z^d} l_k^{(l)}(z)=\sum_{k=1}^{2^l} n_k^{(l)}=n.
\]
We consider now $\C_n^\downarrow(l)$ 
given in \reff{low-eq.11},
and divide it into $\C^{I}(l)$,
where the sum over $i$ spans $\{-1,\dots,-l\}$, 
and $\C^{I\!I}(l)$ for the
remaining terms. In case $j_0>l$, then $\C^{I\!I}(l)$ vanishes. 
Note that for any $h\le L$,
\be{low-BI}
\begin{split}
\sum_{l=h}^{L}
\C^{I\!I}(l)\le &\sum_{l=h}^{L}\sum_{k=1}^{2^{l}}\sum_{i<-l}
2^q\pare{\frac{\alpha_0 \xi_n^\gamma}{2^{l}}}^{q-1}
\pare{ l_{2k}^{(l)}(\D^{(l)}_{2k-1,i})+
l_{2k-1}^{(l)}(\D^{(l)}_{2k,i})}\\
\le & 2^q \sum_{l=h}^{L}\pare{\frac{\alpha_0 \xi_n^\gamma}{2^{l}}}^{q-1}
\sum_{k=1}^{2^l} \sum_{z\in \Z^d} l_{k}^{(l)}(z)\\
\le & n \xi_n \alpha_0^{q-1} \sum_{l\ge 0} \pare{ \frac{1}{2^{q-1}}}^{l}<
\frac{n\xi_n}{16}.
\end{split}
\ee
We have used the condition \reff{alpha0} to obtain the 
last line in \reff{low-BI}.

Now, we use that
\[
P(\sum_{l=h}^{L}\bar\C_n^\downarrow(l)\ge\frac{n\xi_n}{8})\le 
P(\sum_{l=h}^{L}\bar\C^{I}(l)\ge \frac{n\xi_n}{16})+
P(\sum_{l=h}^{L}\bar\C^{I\!I}(l)\ge \frac{n\xi_n}{16}).
\]
Thus, in view of \reff{low-BI}, 
the choice of \reff{alpha0} implies that for any $h$,
$P(\sum_{l=h}^{L}\bar\C^{I\!I}(l)\ge \frac{n\xi_n}{16})=0$.

We proceed now with estimating $\{\bar\C^{I}(l)\}$. 
We do so in three steps. In Step 1, we perform the transformation
of the terms of $\C_n^\downarrow(l)$
into independent variables $X_k$ distributed as geometric variables.
Then, Lemma \ref{lem-5.1}
provides the following bound. For $1/2>\delta>0$, 
\be{bound-key}
P\pare{\sum_{k=1}^{2^{l}} \bar X_k\ge x_{n}} 
\le e^{-\delta x_{n}/4}, \quad\text{if}\quad
4c_u 2^l \max(E[X_k^2]^{1-\delta},E[X_k^2])\le x_{n}. 
\ee
Step 2 establishes the condition in the
right hand side of \reff{bound-key}. Finally, Step 3
compares $x_n$ with the desired rate of decay.

\noindent{\bf Step 1:}
Note that the volume $|\D^{(l)}_{k,i}|$ times the minimal amount of
time spent on sites of $\D^{(l)}_{k,i}$, that is $b_i$,
is bounded by the total time
left for a strand of random walk at generation $l$, so that
\be{low-eq.12}
|\D^{(l)}_{k,i}|\times b_i \le \frac{n}{2^{l}}.
\ee
Now, for fixed $l$ and $0>i\ge-l$, and for $k=1,\dots,2^{l}$, we 
define, following \reff{def-X},
\be{low-eq.14}
X_k:=\pare{b_i\frac{2^{l}}{n}}^{2/d} 
l_{2k-1}^{(l)}(\D^{(l)}_{2k,i})\quad\text{and}\quad
X_k'=\pare{b_i\frac{2^{l}}{n}}^{2/d}
l_{2k}^{(l)}(\D^{(l)}_{2k-1,i}).
\ee
We set $\bar X_k=X_k-E[X_k]$ and $\bar X_k'=X_k'-E[X_k']$,
and using \reff{low-eq.11}, we rewrite
\be{step-amin2}
\acc{\sum_{l=h}^{L}\bar\C^{I}(l)\ge \frac{n\xi_n}{16}}=
\acc{\sum_{l=h}^{L}\sum_{i=-l}^{-1}\sum_{k=1}^{2^{l}} 
2^qb_{i+1}^{q-1}\pare{ \frac{n}{2^{l}b_i}}^{2/d}
(\bar X_k+\bar X_k')\ge \frac{n\xi_n}{16}}.
\ee
We recall an obvious general bound
(for any countable set of indices $\A$) 
\be{step-amin3}
P\pare{ \sum_{a\in \A} X_a\ge \sum_{a\in \A} \alpha_a}\le
 \sum_{a\in \A} P\pare{ X_a\ge \alpha_a}.
\ee
Now, \reff{step-amin3} applied 
to the expression of the right hand side of \reff{step-amin2}, yields
\be{low-eq.13}
P\pare{ \sum_{l=h}^{L}\bar\C^I(l)
\ge \frac{\xi_n n}{16}}\le 2\sum_{l=1}^{L-1}\sum_{i=-l}^{-1}
P\pare{\sum_{k=1}^{2^{l}} \bar X_k\ge x_{n,l,i}},
\ee
with, for $0>i\ge -l$, and $\epsilon>0$ (using 
$\sum_{i\le 0} 2^{\epsilon i}=\sum_{l\ge 0} 2^{-\epsilon l}=
(1-2^{-\epsilon})^{-1}$)
\[
x_{n,l,i}= c_1 \frac{1}{b_{i+1}^{q-1}}
\pare{ \frac{b_i 2^{l}}{n}}^{2/d} 2^{-\epsilon(l-i)}\ n \xi_n,
\quad\text{and}\quad
c_1=\frac{(1-2^{-\epsilon })^2}{32\times 2^{q}}.
\]
The factor 2 appearing in the right hand side of \reff{low-eq.13}
comes from noting that $\bar X_k'$ has the same law, as $\bar X_k$.

Note that for any $k=1,\dots,2^l$, 
$P(X_k>t)\le P(X>t)$, where $X=l_m(\tilde \D)/|\tilde \D|^{2/d}$ with 
$m=n/2^l$, and $\tilde \D$ is a certain level set of
local times $\tilde l_m$ independent from $l_m$.
If $P_0$ and $\tilde P_0$
are the law of the two independent local times, then
Lemma 1.2 of \cite{AC04} yields 
\be{step-amin1}
P(X_k>t)\le
\tilde E_0\cro{ P_0\pare{ l_{m}(\tilde \D)> |\tilde \D|^{2/d} t}}
\le \tilde E_0\cro{ \exp\pare{
-\kappa \frac{ |\tilde \D|^{2/d} t}{|\tilde \D|^{2/d}}}}\le e^{-\kappa t}.
\ee

\noindent{\bf Step 2:}
First, by Lemma~\ref{lem-app.1}, we have
\be{low-eq.22}
E[X_k^2]\le C_d'\psi_d^2(\frac{n}{2^l}) 
\pare{\frac{b_i 2^{l}}{n}}^{4/d}e^{-\kappa_d b_i}.
\ee
When we recall that $b_{i+1}=2b_i$ (see \reff{bitop}),
\reff{bound-key} requires that for some constant $K$, and $0<\delta<1/2$
\be{low-eq.23}
\pare{\psi_d^2(\frac{n}{2^l}) 
\pare{\frac{b_i 2^{l}}{n}}^{4/d}e^{-\kappa_d b_i} }^{1-\delta}
\le K\frac{1}{b_{i}^{q-1}}
\pare{ \frac{b_i 2^{l}}{n}}^{2/d} 2^{-\epsilon(l-i)}\ n \xi_n.
\ee
We use now that $\psi_d^2(k)\le k$ (see \reff{app-key1}),
and \reff{low-eq.23} follows as soon as
\be{low-eq.24}
b_i^{q-1+\frac{2}{d}(1-2\delta)} e^{-(1-\delta)\kappa_d b_i}\le
\pare{ \frac{n}{2^l}}^{\frac{2}{d}(1-2\delta)}2^{-\epsilon(l-i)}\xi_n.
\ee
Recall that $b_i\ge 1$, and the left hand side of \reff{low-eq.24}
is bounded from above and below by constants, when $\delta$ is small
enough. Since $2^l\le
n^{1-\delta_0}$ with $q\delta_0<\frac{2}{d}$, we have that 
\reff{low-eq.24} holds with as soon as $\xi_n\ge 1$
(and choosing $\delta<1/2$ and
$\epsilon$ small enough). In particular, nothing prevents
$\xi_n$ to be as large as possible here.

\noindent{\bf Step 3:} 
Lemma~\ref{lem-bottom} is proved if we show that for some constant $K>0$,
and any $i$ and $l$,
\be{low-eq.20}
\xi_n^{\frac{2}{d}\gamma}n^{1-\frac{2}{d}}\le Kx_{n,l,i}
\quad(\text{recall that }\quad \gamma=\frac{1}{q-1}).
\ee
Condition \reff{low-eq.20} is the most critical to check.
It requires (recalling that $0>i\ge-l$)
\be{low-eq.21}
2^{-i(q-1)} \pare{ \frac{\alpha_0 \xi_n^\gamma 2^{l+i}}
{n}}^{2/d} 2^{-\epsilon (l-i)}n\ge
\xi_n^{\frac{2}{d}\gamma}n^{1-\frac{2}{d}}\Longleftarrow
 (l+i) \frac{2}{d}-i(q-1)>\epsilon(l-i),
\ee
which holds if $2\epsilon< \min(q-1,\frac{2}{d})$.

\epr

\subsection{The {\it top-level} terms}
\label{sec-top}

\bl{lem-top} Assume $d\ge 3$ and $1<q<q_c(d)$. There is a constant $C>0$,
and $\delta>0$,
such that for any $h\in \{1,\dots,L\}$ and $\xi_n\ge n^{-\delta}$
\be{top.1}
P\pare{ \sum_{l=h}^L \bar \C_n^\uparrow(l)\ge \frac{\xi_n n}{8}}\le
\exp\pare{ -C \xi_n^{\frac{2}{d}\frac{1}{q-1}}
\min(1,\xi_n^{2/d}) n^{1-\frac{2}{d}}}.
\ee
\el

\bpr
Following \cite{A06}, we take two sequences of
positive numbers $\{a_i,\ i=1,\dots,M_n\}$, and for each $i$
$\{p_l^{(i)},\ l=h,\dots,L-1\}$ (to be made explicit later) with
\be{bulk.4}
\sum_{i=1}^{M_n} a_i=1,\quad\text{and for each $i$}\quad
\sum_{l=h}^{L-1} p_l^{(i)}=1.
\ee
For any $h\le L$, we have
\be{top.2}
\begin{split}
P\big(\sum_{l=h}^{L}&\bar\C_n^\uparrow(l)
\ge\frac{\xi_n n}{8}\big)\\
&\le 2\sum_{l=1}^{L}\sum_{i\ge 1}
P\pare{\sum_{k=1}^{2^{l}} \ind_{\G^{(l)}_{2k,i}}
l_{2k-1}^{(l)}(\D^{(l)}_{2k,i})-
E\cro{ \ind_{\G^{(l)}_{2k,i}}l_{2k}^{(l)}(\D^{(l)}_{2k-1,i})}
\ge\frac{n}{8\ 2^q \beta_{i+1}^{q-1}}p_l^{(i)}a_i}.
\end{split}
\ee
We proceed as in the proof of Lemma
\ref{lem-bottom} with Steps 1,2 and 3.

\noindent{\bf Step 1:}
We first bound $|\D_{2k,i}^{(l)}|$.
Note that on $\G^{(l)}_{2k,i}$ for any $k,i$, we have 
\be{bulk.5}
|\D_{2k,i}^{(l)}|\le \min\pare{ \frac{n/2^l}{\beta_i \xi_n^\gamma},
\frac{n(2\kappa(q,d)+\xi_n)}{\beta_i^q\xi_n^{\gamma+1}}}\le
\frac{n}{\beta_i \xi_n^\gamma\min(1,\xi_n)}\min\pare{\frac{1}{2^l}, 
\frac{2\kappa(q,d)+1}{\beta_i^{q-1}}}.
\ee
We used in \reff{bulk.5} that
\[
\frac{2\kappa(q,d)+\xi_n}{\xi_n} \le 
(2\kappa(q,d)+1)\frac{\max(\xi_n,1)}{\xi_n}=
\frac{2\kappa(q,d)+1}{\min(\xi_n,1)}.
\]
In order to use Lemma \ref{lem-5.1}, we need to normalize
$l^{(l)}_{2k}\pare{ \D_{2k-1,i}^{(l)}}$ with a constant
smaller than $|\D_{2k-1,i}^{(l)}|^{-2/d}$. We choose, for $l$ and $i$ fixed,
\be{bulk.8}
\zeta^{(l)}_{i}=\pare{\frac{\beta_i \xi_n^\gamma\min(1,\xi_n)}
{n}}^{\frac{2}{d}}
\left\{ \begin{array}{ll}
                (2\kappa(q,d)+1)^{-\frac{2}{d}}\beta_i^{\frac{2}{d}(q-1)}&
\mbox{ for } l\le l^*_i \, ,
                \\
                2^{\frac{1}{d}l} & \mbox{ for } l> l^*_i \, ,
                \end{array}
        \right.
\ee
with $l^*_i$ is such that $2^{l^*_i}=(2\kappa(q,d)+1)^{-2}\beta_i^{2(q-1)}$.
As in \reff{step-amin1}, we set
\[
X_k=\zeta^{(l)}_{i} \ind_{\acc{\G_{2k,i}^{(l)}}}
l^{(l)}_{2k-1}\pare{ \D_{2k,i}^{(l)}}, \quad\text{and}\quad
P(X_k>t)\le \exp(-\kappa t).
\]
Using \reff{bulk.5},
and the notation $\bar X_k$ for $X_k-E[X_k]$, we have
\be{bulk.6}
P\pare{ \sum_{l=h}^{L}\bar\C_n^\uparrow(l)
\ge \frac{\xi_n n}{8}}\le 2
\sum_{l=h}^{L}\sum_{i\ge 1} P\pare{\sum_{k=1}^{2^{l}}\bar X_k\ge x_{n,l,i}}
\text{  with  }
x_{n,l,i}=\frac{n\zeta^{(l)}_{i}}{16(2^q+1)\beta_{i+1}^{q-1}}a_i p_l^{(i)}.
\ee
\noindent{\bf Step 2:} We establish now a condition equivalent to
\reff{bound-key}.
\be{i-ii}
2^{l(1+\delta_0\frac{2}{d})} E[X_k^2]\le K 
\xi_n^{\frac{2}{d}\gamma}\min(1,\xi_n^{2/d}) n^{1-\frac{2}{d}},
\ee
for some
constant $K$. Thus, when $l\le l_i^*$, and for some constant $C$
\be{bulk-bis}
\begin{split}
2^{l(1+\delta_0\frac{2}{d})} E[X_k^2]\le&2^l2^{l\delta_0\frac{2}{d}}
\pare{\beta_i^q\frac{\xi_n^\gamma}{n}}^{4/d}
\min(1,\xi_n^{4/d})C_d\psi_d^2(\frac{n}{2^l})
e^{-\kappa_d \xi_n^\gamma \beta_i} \\
\le & C_d \min(1,\xi_n^{4/d})
\frac{n^{1-\frac{2}{d}}}{\xi_n^{4/d}}
\pare{\frac{2^l}{n} \psi_d^2(\frac{n}{2^l})}
\pare{ \frac{2^{l\delta_0}}{n}}^{\frac{2}{d}}
\sup_{x>0}\acc{x^{\frac{4q}{d}} \exp(-\kappa_d  x)}\\
\le & C n^{1-\frac{2}{d}-\frac{2}{d}(1-\delta_0(1-\delta_0))}.
\end{split}
\ee
In this case, \reff{i-ii} holds if for some constant $C$
\be{small-xi.1}
\xi_n^{\gamma}\min(1,\xi_n)\ge \frac{C}{n^{(1-\delta_0(1-\delta_0))}}.
\ee
\reff{small-xi.1} holds when $\xi_n\ge n^{-\delta}$ for $\delta>0$
small enough. When $l>l_i^*$, for a constant $C'$
\be{bulk.14}
\begin{split}
2^{l(1+\delta_0\frac{2}{d})} E[X_k^2]\le&2^l2^{l\delta_0\frac{2}{d}} 
\pare{\beta_i\frac{\xi_n^\gamma}{n}}^{4/d}
\min(1,\xi_n^{4/d})C_d\psi_d^2(\frac{n}{2^l})
e^{-\kappa_d \xi_n^\gamma \beta_i} 2^{\frac{2}{d}l}\\
\le & C_d n^{1-\frac{2}{d}}\min(1,\xi_n^{4/d})
\pare{\frac{2^l}{n} \psi_d^2(\frac{n}{2^l})}
\pare{ \frac{2^{l(1+\delta_0)}}{n}}^{\frac{2}{d}}
\sup_{x>0}\acc{x^{\frac{4}{d}} \exp(-\kappa_d  x)}\\
\le & C' \min(1,\xi_n^{4/d}) n^{1-\frac{2}{d}-\frac{2}{d}\delta_0^2}.
\end{split}
\ee
When $\xi_n\ge 1$, \reff{i-ii} follows from \reff{bulk-bis}
and \reff{bulk.14}. When $\xi_n< 1$, we need in addition that
$\xi_n^{\gamma-1}\ge n^{-\delta_0^2}$.

\noindent{\bf Step 3:}
We show that we can choose $p_l^{(i)}$ and $a_i$ such that 
for any $i,l$ (and $n$ large enough), there is
a constant $c$, independent on $i,l$ and $n$, and
\be{bulk.10}
\frac{n\zeta^{(l)}_{i}}{\beta_{i+1}^{q-1}}
p_l^{(i)}a_i \ge cn^{1-\frac{2}{d}}\xi_n^{\frac{2}{d}\gamma}
\min(1,\xi_n^{2/d}).
\ee
It is possible to choose
a normalizing constant $a_0$ (which depends only on $q$), 
such that for $i=1,\dots,M_n$, $\sum_{i\ge 1} a_i\le 1$, where
\be{bulk.11}
a_i=a_0\pare{ \frac{\beta_{i+1}^{q-1}}{\beta_i^{\frac{2}{d}q}} }^{1/2}=
a_0 2^{\frac{q-1}{2}}
2^{ -\alpha i},\quad\text{with}\quad
\alpha:=\frac{1}{2}\pare{1-\frac{q}{q_c(d)}}. 
\ee
Indeed, the condition $q<q_c(d)$ implies that $\alpha$ is positive,
and the series in \reff{bulk.11} is convergent.

Now, for a fixed $i\ge 1$, we turn to the choice of $\{p_l^{(i)},l\ge 1\}$.
We will choose later two constants $p^*$ and $\bar p$
such that for $l<l^*_i$,
\be{bulk.12}
p_l^{(i)}=p^* 2^{ -\alpha i},
\ee
whereas for $l>l_i^*$, 
\be{bulk.13}
\begin{split}
p_l^{(i)}=&\bar p \frac{\beta_{i+1}^{q-1}}{\beta_i^{2/d}} 
\frac{2^{\alpha i}}{2^{l/d}}
\le \bar p \frac{\beta_{i+1}^{q-1}}{\beta_i^{2/d}}\pare{ 
\frac{\beta_{i}^{\frac{2}{d}q}}{\beta_{i+1}^{(q-1)}} 
\frac{1}{2^{l_i^*/d}}}^{1/2} 
\frac{2^{\frac{q-1}{2}}}{2^{l/(2d)}}\\
\le & \bar p \pare{ \frac{\beta_{i+1}^{q-1}}{\beta_i^{\frac{2}{d}q}}
\frac{2^{q-1}}{2^{(l-l_i^*)/d}}}^{1/2}
\le  \bar p  \frac{2^{\frac{q-1}{2}}}{2^{(l-l_i^*)/(2d)}}.
\end{split}
\ee
We proceed now to normalize $\{p_l^{(i)},l\ge 1\}$. We need
to choose $p^*$ and $\bar p$ such that for each $i$,
$\sum_l p_l^{(i)}\le 1$. Recall that there is $c_1$ such that $l_i^*\le
c_1i$. Now, note that
\ba{no-weak17}
\sum_{l} p_l^{(i)}&\le&p^* l^*_i 2^{-\alpha i}+
\bar p 2^{\frac{q-1}{2}}\sum_{l>l_i^*}  2^{-(l-l_i^*)/(2d)}\cr
&\le & p^* c_1i 2^{-\alpha i}+
\bar p 2^{\frac{q-1}{2}}\sum_{l>0} \frac{1}{2^{l/(2d)}}\cr
&\le & c_1 p^* \sup_{x>0}\acc{x 2^{-\alpha x}}
+\bar p\frac{2^{\frac{q-1}{2}}}{2^{1/(2d)}-1}.
\ea
It is important to see that $p^*$ and $\bar p$ can be chosen
independently of $i$.
Now, we check \reff{bulk.10}. For $l<l_i^*$,
\be{top.3}
\frac{n\zeta^{(l)}_{i}}{\beta_{i+1}^{q-1}} p_l^{(i)}a_i 
=a_0 p^* \frac{n\zeta^{(l)}_{i}}{\beta_i^{\frac{2q}{d}}}
=a_0 p^* \xi_n^{\frac{2}{d}\gamma}\min(1,\xi_n^{2/d})
 n^{1-\frac{2}{d}}.
\ee
For $l\le l_i^*$,
\be{top.4}
\frac{n\zeta^{(l)}_{i}}{\beta_{i+1}^{q-1}} p_l^{(i)}a_i
=a_0 \bar p 2^{(q-1)/2} \frac{n\zeta^{(l)}_{i}}
{2^{l/d}\beta_i^{\frac{2}{d}}}=a_0 \bar p 2^{(q-1)/2} 
\xi_n^{\frac{2}{d}\gamma} \min(1,\xi_n^{2/d}) n^{1-\frac{2}{d}}.
\ee
This concludes the proof of Lemma~\ref{lem-top}.

\epr

\subsection{The lower bound in \reff{intro.4}}\label{sec-LBd3}
As in inequalities (80) and (81) of~\cite{AC05},
the lower bound follows from H\"older's inequality.
Indeed, it is immediate that $||l_n||_q^q/n\ge (n/|R_n|)^{q-1}$, where $R_n$ is
the set of visited sites up to time $n$. Thus, when $n$ is large enough
\[
|R_n|\le \frac{n}{(2\kappa(q,d)+\xi_n)^{\gamma}}
\Longrightarrow 
||l_n||_q^q\ge n(2\kappa(q,d)+\xi_n)\ge E[||l_n||_q^q]+\xi_n n.
\]
Now, forcing the walk to stay in a ball $B(0,r_n)$
centered at the origin, and of radius $r_n$ satisfying 
$r_n^d=n/(2\kappa(q,d)+\xi_n)^{\gamma}$ 
implies that $|R_n|\le n/(2\kappa(q,d)+\xi_n)^{\gamma}$. The cost of this constraint is 
$\exp(-c \frac{n}{r_n^2})$, which yields the lower bound in \reff{intro.4}, when
we recall that $\xi_n\ge 1$.

\subsection{Proof of Theorem~\ref{theo-d3}}\label{sec-proofd3}
We collect the estimates of the previous subsections in order
to prove \reff{intro.4} allowing $\xi$ to depend on $n$, as in
Remark~\ref{rem-xin}.

When $\xi_n\ge 1$, using the decomposition
\reff{bulk.3}, with the estimates \reff{cor-LD}, \reff{bottom.1}
and \reff{top.1}, we obtain the upper bound in \reff{intro.5}.
The lower bound in \reff{intro.5} follows from Section~\ref{sec-LBd3}.

When $\xi_n<1$, then Lemma~\ref{lem-LD} imposes that $\xi_n\ge n^{-\delta}$
with $0\le \delta<\delta_0/2$, whereas Lemma~\ref{lem-top} holds
for some positive $\delta$. Thus, we conclude that Remark~\ref{rem-xin}
holds with \reff{eq-xin1}. Note that a lower bound is missing in this
case.
\section{The super-critical regime}\label{sec-supercritical}
We consider here $q>q_c(d)=\frac{d}{d-2}$.
The main result of this section is to show that
only sites of $\{z:l_n(z)\ge(n\xi_n)^{1/q}/A\}$ (for some $A>0$),
contribute to realize the excess self-intersection, at a cost
given in \reff{intro.6}.

The proof of Theorem~\ref{theo-d4} relies on the
following estimates. For any $\epsilon$ with 
$0<\epsilon<1/q$, and any $\delta$, with $0<\delta<1/3$, 
and two constants $A>A_0$, we write
\be{super-b1}
\begin{split}
P\big( ||l_n||_q^q-E[||l_n||_q^q]\ge&\xi_n n\big)\le
P\pare{ \sum_{z} \ind\acc{ z:\ l_n(z)<\xi_n^{1/q} n^{1/q-\epsilon}}
l_n^q(z) -E[||l_n||_q^q]\ge n\delta\xi_n}\\
&+ P\pare{ \sum_{z} \ind\acc{\xi_n^{1/q} n^{1/q-\epsilon}<l_n(z)
\le \frac{(\xi_n n)^{1/q}}{A}}\ l_n^q(z) \ge n\delta\xi_n}\\
&+P\pare{ \sum_{z} \ind\acc{\frac{(\xi_n n)^{1/q}}{A}<l_n(z)
\le \frac{(\xi_n n)^{1/q}}{A_0}}\ l_n^q(z)\ge n\xi_n(1-3\delta)}\\
&+P\pare{ \sum_{z} \ind\acc{l_n(z)
> \frac{(\xi_n n)^{1/q}}{A_0}}\ l_n^q(z)\ge n\xi_n\delta}.
\end{split}
\ee
In Section~\ref{sec-superlow}, we show that the contribution
of $\{z:\ l_n(z)<\xi_n^{1/q}n^{1/q-\epsilon}\}$, for any $\epsilon>0$,
is negligible. More precisely, we establish that 
there is $\epsilon'>0$ such that for any $\delta>0$, and $n$
large enough
\be{super-step1}
P\pare{ \sum_{z\in \Z^d}
\ind\acc{ z:\ l_n(z)<\xi_n^{1/q} n^{1/q-\epsilon}}
l_n^q(z) -E[||l_n||_q^q]\ge n\delta\xi_n}\le
\exp\pare{ -\xi_n^{1/q} n^{\frac{1}{q}+\epsilon'}}.
\ee
The proof of \reff{super-step1}
is similar to that of Theorem \ref{theo-d3}.

In Section~\ref{sec-supermid}, we show the following lemma.
\bl{lem-topII} Assume $d\ge 3$ and $q>q_c(d)$. There is
constants $A_0$ and $\kappa_d$, such that for $\epsilon>0$, and
any $\xi_n>0$, and any integer $n$,
\be{topII-main}
P\pare{ \sum_{z\in \Z^d} \ind\acc{\xi_n^{1/q} n^{1/q-\epsilon}<l_n(z)
\le \frac{(\xi_n n)^{1/q}}{A_0}}\ l_n(z)^q\ge n\xi_n}\le
\exp\pare{ -\kappa_d  \xi_n^{\frac{1}{q}}\ n^{\frac{1}{q}}}.
\ee
Furthermore, there is a constant $C>0$ such that
for $\delta>0$, and $A>A_0$ 
\be{super-eq30}
P\pare{ \sum_{z\in \Z^d} \ind\acc{\xi_n^{1/q} n^{1/q-\epsilon}<l_n(z)\le
\frac{(\xi_n n)^{1/q}}{A}}\ l_n^q(z) \ge n\delta\xi_n}\le
\exp\pare{ -C A \delta^{1-\frac{2}{d}}n^{\frac{1}{q}}}.
\ee
\el
Finally, since we have a transient random walk,
it is obvious that for $c>0$,
\[
P\pare{ \sum_{z} \ind\acc{l_n(z)
\ge \frac{(\xi_n n)^{1/q}}{A_0}}\ l_n^q(z)\ge n\xi_n\delta}\le
P\pare{\exists z: l_n(z) \ge \frac{(\xi_n n)^{1/q}}{A_0}}\le
ne^{-c\frac{(\xi_n n)^{1/q}}{A_0}}.
\]

The lower bound comes from requiring that the origin
is visited $(n \xi_n)^{1/q}$ times.

\subsection{The contribution of $\{z:\ l_n(z)<\xi_n^{1/q} n^{1/q-\epsilon}\}$}
\label{sec-superlow}
The first step is to perform a approximation of $||l_n||_q^q$
over $\{z:\ l_n(z)<\xi_n^{1/q} n^{1/q-\epsilon}\}$ as in
Section \ref{sec-general}. This is explained in Remark \ref{rem-legall2}.

To allow for the possibility of $\xi_n$ to depend on $n$, we need
to trace the occurrences of $\xi_n$, and in this respect, it is useful to
modify the subdivision chosen in \reff{bitop}. 
We choose again
$\alpha_0$ as in \reff{alpha0}, and for $i\ge -j_0$ we keep 
$\beta_i=\alpha_0 2^i$, and
\be{bibottom}
\forall i<0\quad b_i=\xi_n^{\frac{1}{q-1}}\beta_i
\quad\text{and}\quad \forall i\ge 0\quad b_i=\xi_n^{\frac{1}{q}}\beta_i.
\ee
Recall that when $\xi_n<1$, then $\D^{(l)}_{k,i}=\emptyset$ 
for $i<0$, and $\C_n^\downarrow$ vanishes.
However, when $\xi_n\ge 1$, for each $k$ and $l$,
there is an overlap between $\D^{(l)}_{k,-1}$
and $\D^{(l)}_{k,0}$ since $\xi_n^{1/q}\le \xi_n^{1/(q-1)}$.

For a small $\epsilon>0$, the subdivision $\{b_i\}$
covers $[1,\xi_n^{1/q} n^{1/q-\epsilon}]$.
As in the proof of Theorem~\ref{theo-d3}, we start with \reff{bulk.3}.
We first treat $\C_n^\uparrow(l)$.
\bl{lem-topI} Assume $d\ge 3$, and $q>q_c(d)$. We consider
a sequence $\{\xi_n,n\in \N\}$ such that for some
$\delta>0$ small enough $\xi_n\ge n^{-\delta}$. 
There is a constant $\epsilon'>0$,
such that for any $h\in \{1,\dots,L\}$ and for $n$ large enough
\be{super-eq1}
P\pare{ \sum_{l=h}^L \bar \C_n^\uparrow(l)\ge \frac{\xi_n n}{8}}\le
\exp\pare{ -\xi_n^{\frac{1}{q}}\min(1,\xi_n^{2/d}) 
n^{\frac{1}{q}+\epsilon'}}.
\ee
When $q=q_c(d)$, then for any $h\in \{1,\dots,L\}$, and 
$n$ large enough
\be{super-crit}
P\pare{ \sum_{l=h}^L \bar \C_n^\uparrow(l)\ge \frac{\xi_n n}{8}}\le
\exp\pare{ -\xi_n^{\frac{1}{q}}\min(1,\xi_n^{2/d}) 
n^{\frac{1}{q}-\epsilon'}}.
\ee

\el
\br{rem-superxi-n}
When $1> \xi_n\ge n^{-\delta}$ with $\delta$ small, then
the terms $\{\C_n^\downarrow(l),l\le L\}$ vanish, whereas 
$S_q^{(L)}$ is negligible. Indeed, according to Lemma~\ref{lem-LD},
it suffices to show that
\[
\xi_n n^{1-q\delta_0}\ge \xi_n^{1/q+2/d} n^{1/q+\epsilon'},
\]
which holds when $\xi_n> 1/\sqrt{n}$ (which we always assume). 

When $\xi_n\ge 1$, and for a choice of $\delta_0<2/(dq)$, we have
$\xi_n n^{1-q\delta_0}\ge (\xi_n n)^{1/q}$ so that by
\reff{self-key1}, we can neglect $S_q^{(L)}$. Also,
recall that we can assume $\xi_n\le n^{q-1}$ (see Remark~\ref{rem-growth}).
This latter bound is equivalent to 
\[
\xi_n^{\frac{2}{d}\frac{1}{q-1}} n^{1-\frac{2}{d}}\ge (\xi_n n)^{1/q}.
\]
Now, \reff{bottom.1} of Lemma~\ref{lem-bottom}
allows us to neglect the contribution of $\{\C_n^\downarrow(l),l\le L\}$.
\er

\noindent{{\bf Proof of Lemma~\ref{lem-topI}}}
We proceed with Steps 1,2 and 3 as in the proofs of
Lemma~\ref{lem-bottom} and Lemma~\ref{lem-top}.

\noindent{\bf Step 1:} 
The first difference with the proof of Theorem \ref{theo-d3},
is the choice of the subdivision $\{b_i\}$ of \reff{bibottom}.
Note that the bound on $|\D_{k,i}^{(l)}|$ of \reff{bulk.5} becomes
\[
|\D_{k,i}^{(l)}|\le (2\kappa(q,d)+1) \frac{n}{\min(1,\xi_n)\beta_i^q}.
\]
This implies a new definition for $\zeta_i^{(l)}$.
Also, note that the choice \reff{bulk.11} for $a_i$ is not
possible since $\alpha<0$ in this case. Thus, we set for $i\in \N$,
and $\delta>0$ to be chosen later,
\be{super-eq2}
a_i=(1-2^{-\delta}) 2^{-\delta i},\quad
p_l^{(i)}=(1-2^{-\delta}) 2^{-\delta l},\quad
\text{and}\quad
\zeta_i^{(l)}=\pare{\frac{\beta_i^q}{n}\min(1,\xi_n)}^{2/d}.
\ee
Accordingly, inequality \reff{bulk.6} holds, but with
\be{new-xn}
\begin{split}
x_{n,l,i}=\frac{n\xi_n^{1/q}\min(1,\xi_n^{2/d})
\zeta^{(l)}_{i}}{16(2^q+1)\beta_{i+1}^{q-1}} p_l^{(i)}a_i
=&c 2^{-\delta (i+l)} \beta_i^{q\frac{2}{d}-(q-1)}
\xi_n^{1/q}\min(1,\xi_n^{2/d}) n^{1-\frac{2}{d}}\\
=&c2^{-\delta (i+l)}\xi_n^{1/q}\min(1,\xi_n^{2/d}) 
n^{1/q_c(d)}\beta_i^{1-q/q_c(d)} .
\end{split}
\ee
Note that $q>q_c(d)$ implies that $x_{n,l,i}$ is small 
when $\beta_i$ is large.

\noindent{\bf Step 2:} 
We establish that
for $\delta>0$ small $2^{l(1+\delta)}E[X_k^2]\le x_{n,l,i}$. 
This latter inequality is equivalent to 
\be{step-2}
2^{l(1+\delta)} \zeta^{(l)}_{i} \psi_d^2(\frac{n}{2^l})
e^{-\kappa_d \xi_n^{1/q}\beta_i}\le c \frac{n\xi_n^{1/q}}{\beta_i^{q-1}}
2^{-\delta(i+l)},
\ee
which is equivalent to
\be{new-i}
\pare{ \frac{2^l}{n} \psi_d^2(\frac{n}{2^l})}2^{\delta (i+2l)}
\beta_i^{q\frac{2}{d}+(q-1)} e^{-\kappa_d \xi_n^{1/q}\beta_i}\le
c n^{2/d} \xi_n^{1/q}\min(1,\xi_n^{2/d}).
\ee
Since $\psi_d^2(k)\le k$ when $d\ge 3$, \reff{new-i} holds for any
$\beta_i$, $\delta>0$ small enough, and $\xi_n\ge n^{-\delta}$. 

\noindent{\bf Step 3:}
We distinguish the
cases $q>q_c(d)$ and $q=q_c(d)$. 

When $q>q_c(d)$, then we need to show that
\[
x_{n,l,i}\ge \xi_n^{1/q} \min(1,\xi_n^{2/d})n^{1/q+\epsilon'}.
\]
using \reff{new-xn}, this is equivalent to
\be{super-eq5}
n^{1/q_c(d)}\beta_i^{1-q/q_c(d)} 2^{-\delta (i+l)}
\ge c n^{\frac{1}{q}+\epsilon'}.
\ee
So \reff{super-eq5} holds if for some $\epsilon'>0$
\be{super-eq7}
2^{\delta (i+l)}\beta_i^{q/q_c(d)-1}\le n^{1/q_c(d)-1/q-\epsilon'}.
\ee
Since $\beta_i\le n^{1/q-\epsilon}$, \reff{super-eq7} holds
for $\delta$ and $\epsilon'$ both small enough.

When $q=q_c(d)$, we need to show that
\[
x_{n,l,i}\ge\xi_n^{1/q} \min(1,\xi_n^{2/d}) n^{1/q_c(d)-\epsilon'}.
\]
This is obvious
as soon as $n^{\epsilon'}\ge 2^{\delta(l+i)}$, which holds for
$\epsilon'>0$, when $\delta$ is small enough.

\qed
\subsection{Proof of Lemma~\ref{lem-level}}\label{sec-gamma}
Since the proof of Lemma~\ref{lem-level} is similar to the
proof given in Section \ref{sec-superlow}, we do not give all
details, but only focus on the differences. 
When dealing with $\{||\Theta_{n^b}(l_n)||_q^q\ge n^a\}$,
with $a>1$, our starting point is inequality \reff{bootstrap.2} of
Corollary~\ref{cor-legall} with $M=n^b$. We choose $m_n=E[||l_n||_q^q]
+\epsilon n^a$, for $\epsilon<1/2$,
and $\epsilon_n=(1-\epsilon) n^a$. 
We use the sets $\{\D_{k,i}^{(l)},\ i\in \N\}$ 
of Section \ref{sec-superlow}. 
However, $\{b_i,i\in \N\}$ only cover
$[1,n^b]$, and $\xi_n$ of \reff{bitop} is set to 1. This latter
choice implies that there is no term $\C_n^\downarrow(l)$. 
Note that the bootstrap bound of \reff{bulk.1}
defining $\G_{k,i}^{(l)}$ is here $\{|\D_{k,i}^{(l)}|\le n^a/\beta_i^q\}$. 

Now, we proceed as in the proof of Lemma~\ref{lem-level}.
To see that $S_q^{(L)}$ has a negligible
contribution, note that for $a>1$ and any $\epsilon>0$, 
\reff{self-key1} implies that
\[
P\pare{S_q^{(L)}(n)-E[||l_n||_q^q]\ge \epsilon n^a}\le
e^{-\epsilon n^{a-q\delta_0}}.
\]
Since $q \delta_0<2/d$, it is enough (and easy) to check that
for $a>1$ and $q\ge q_c(d)$
\[
a-\frac{2}{d}>(1-\frac{2}{d})a-(\frac{q}{q_c(d)}-1)b.
\]
The main differences with the proof of Section \ref{sec-superlow},
is $\zeta^{(l)}_{i}$ and $x_{n,l,i}$ which read here
\be{zeta-gamma}
\zeta^{(l)}_{i}=\frac{\beta_i^{\frac{2}{d}q}}{n^{\frac{2}{d}a}},
\quad\text{  and  }\quad x_{n,l,i}=\frac{n^a\zeta^{(l)}_{i}}{2(2^q+1)
\beta_{i+1}^{q-1}} p_l^{(i)}a_i.
\ee
Step 2 (similar to \reff{step-2})
is easy to check here, and we omit to do it. 

To check Step 3, i.e. the condition corresponding to \reff{low-eq.20}, 
we recall the definition of $a_i$ and $p_l^{(i)}$ given in \reff{super-eq2},
and use that $\beta_i\le n^b$ (and $q\ge q_c(d)$), to obtain
\be{xn-gamma}
x_{n,l,i}=c2^{-\delta (i+l)} n^{a(1-\frac{2}{d})}\beta_i^{1-q/q_c(d)}
\ge c 2^{-\delta (i+l)} n^{a(1-\frac{2}{d})-b(q/q_c(d)-1)}.
\ee
In conclusion, we obtain for any $\epsilon>0$, and $\delta>0$ small enough
\be{super-eq3}
P\pare{ \sum_{l=h}^L \bar \C_n^\uparrow(l)\ge  n^a}\le
\exp\pare{-n^{\zeta(q,a,b)-\epsilon}}\text{   with   }
\zeta(q,a,b)=a(1-\frac{2}{d})-b(\frac{q}{q_c(d)}-1).
\ee

\subsection{The contribution of $\{z:\ \xi_n^{1/q}
n^{1/q-\epsilon} < l_n(z)\}$} \label{sec-supermid}
In this section, we prove Lemma~\ref{lem-topII}.
We deal with sites whose local times is close to
$n^{1/q}$. 
We follow now the proof of Lemma 3.1 of \cite{A07}. 
Let $\{\alpha_i,i=1,\dots,M\}$ be a subdivision
of $[\frac{1}{q}-\epsilon,\frac{1}{q}]$, to be chosen later.
We justify later in the proof, the choice of
\be{new-A0}
A_0=\exp\pare{2\pare{ \sqrt{ \frac{q}{q_c(d)}}-1}}.
\ee
Also, let $\{p_i,i=0,\dots,M\}$ be positive number
summing up to 1, and define for $i<M$, and $A\ge A_0$
\be{super-eq38}
\D_i=\acc{z: \xi_n^{1/q}n^{\alpha_i}\le l_n(z)< 
\xi_n^{1/q}n^{\alpha_{i+1}}},\quad\text{and}
\quad \alpha_M=\frac{1}{q}-\frac{\log(A)}{\log(n)}.
\ee
Now, as in (3.5) of \cite{A07} (see also Lemma 3.1 of \cite{AC05}), 
we have for any $\delta>0$
\be{super-eq39}
P\pare{ \sum_{z\in\cup \D_i} l_n^q(z)\ge n\delta\xi_n}\le \sup_{0\le i< M}
\acc{ C_i(n) \exp\pare{ -\kappa_d \xi_n^{1/q}\delta^{1-\frac{2}{d}}
n^{\zeta_{i}}p_i^{1-\frac{2}{d}}}},
\ee
with an innocuous combinatorial term $C_i(n)$ independent on $\xi_n$. 
For $0\le i<M$,
\be{super-eq40}
\begin{split}
\zeta_{i}&=\alpha_i+(1-\frac{2}{d})(1-q\alpha_{i+1})\\
&=\frac{1}{q}+ \frac{q}{q_c}\pare{ \frac{1}{q}-\alpha_{i+1}}-\pare{ \frac{1}{q}-\alpha_{i}}.
\end{split}
\ee
Set $a=\sqrt{q/q_c}>1$, and for $i<M$
\be{super-eq41}
\frac{1}{q}-\alpha_{i}=a\pare{ \frac{1}{q}-\alpha_{i+1}},
\quad\text{so that}\quad
\frac{1}{q}-\alpha_{i}=a^{M-i}\pare{ \frac{1}{q}-\alpha_{M}}=
\frac{a^{M-i}\log(A)}{\log(n)}.
\ee
Now, $M$ is chosen such that $\alpha_0=\frac{1}{q}-\epsilon$, that is
$\epsilon\log(n)=a^M\log(A)$. Also, we have
$\zeta_i=\frac{1}{q}+ (a-1) \pare{ \frac{1}{q}-\alpha_{i}}$, and 
we choose (with a normalizing constant $\bar p$ depending only on $a$)
\be{super-eq42}
\pare{\frac{p_i}{\bar p}}^{1-\frac{2}{d}}
=e^{-(a-1)a^{M-i}}\quad\text{and}\quad
n^{\zeta_i}p_i^{1-\frac{2}{d}}=n^{\frac{1}{q}}\bar p^{1-\frac{2}{d}}
e^{\pare{\pare{\log(A)-(a-1)}a^{M-i}}}.
\ee
We need to choose $\log(A_0)>(a-1)$, and our arbitrary choice
of \reff{new-A0} achieves this goal. Thus, the smallest value of
$n^{\zeta_i} p_i^{1-\frac{2}{d}}$ is $n^{\frac{1}{q}}\bar p A\exp(1-a)$.
When we choose $A=A_0$, and $\delta=1$,
we obtain \reff{topII-main}, whereas when
we choose $A>A_0$, and $\delta<1$, we reach \reff{super-eq30}.

\section{About the CLT.}\label{sec-clt}
It will be convenient to use, in this section, the notation $\L_q(n)=
||l_n||_q^q$.
\subsection{Expectation Estimates.}
\noindent{{\bf Proof of Lemma~\ref{lem-moy}}}

Let $n_1$ and $n$ be two integers with $n_1\le n$, and let 
$n_2=n-n_1$. Taking expectation in \reff{decomp.3} yields
\be{moy.2}
E[S_q^{(1)}]\le E[\L_q(n)]\le E[S_q^{(1)}]+E[\I_1(n_1,n_2)].
\ee
We choose a subdivision $\{b_i,i\in \N\}$ with $b_i=i$,
and compute $E[\I_1(n_1,n_2)]$.
Now, using inequality \reff{app-key2} of Lemma~\ref{lem-app.1},
as well as \reff{quasi-dyadic} we have constants $c_d$ such that,
when calling $l_{n_1}^{(1)}=l^{(1,1)}_{]0,n_1]}$ and
$l_{n_2}^{(2)}=l^{(1,2)}_{[0,n_2[}$, 
and using that the local time of a site
increases with the length of the time-period,
\be{var.3}
\begin{split}
E[\I_1(n_1,n_2)]\le & 2^{q}\sum_{z\in \Z^d}\sum_{i\ge 1} b_{i+1}^{q-1}
\pare{l_{n_1}^{(1)}\pare{\acc{z:l_{n_2}^{(2)}(z)\ge b_i}}+
l_{n_2}^{(2)}\pare{\acc{z:l_{n_1}^{(1)}(z)\ge b_i}}}\\
\le& C_d\ \psi_d(\max(n_1,n_2))\sum_{i\ge 1} (i+1)^{q-1} e^{-\kappa_d i}
\le c_d\ \psi_d(\max(n_1,n_2)).
\end{split}
\ee
Thus, if we call $a(n)=E[\L_q(n)]$, and use \reff{moy.2} and \reff{var.3}
\be{var.21}
a(n_1)+a(n_2)\le a(n)\le a(n_1)+a(n_2)+ c_d\psi_d(\max(n_1,n_2)).
\ee
We fix an integer $n$, and for any $k$ (going to infinity), we perform
its euclidean division $k=m_kn+r_k$ with $0\le r_k<n$, and obtain from \reff{var.21}
\be{var.22}
m_k a(n)\le m_k a(n)+a(r_k)\le a(m_kn+r_k)\le a(m_kn)+a(r_k)+c_d \psi_d(m_kn).
\ee
Now, we can use the almost dyadic decomposition of $m_k$, so that if
$L(m_k)$ denote the integer part of $\log_2(m_k)+1$, we have
\be{var.23}
\begin{split}
a(m_kn)\le&a(m_1^{(1)}n)+a(m_2^{(1)}n)+c_d(\psi_d(m_1^{(1)}n)+\psi_d(m_2^{(1)}n))\\
\le & m_k a(n)+c_d \sum_{l=1}^{L(m_k)}\sum_{j=1}^{2^l} \psi_d(m_j^{(l)}n)\\
\le & m_k a(n)+2c_d \sum_{l=1}^{L(m_k)} 2^l \psi_d\pare{\frac{m_k}{2^l}n}\\
\le & m_k a(n)+4c_d \psi_d(n) m_k.
\end{split}
\ee
The last line of \reff{var.23} is obtained after a
simple computation that we omit.
Thus, we are left with
\be{var.24}
\frac{nm_k}{nm_k+r_k}\ \frac{a(n)}{n}\le
\frac{a(k)}{k}\le \frac{nm_k}{nm_k+r_k}\  \frac{a(n)}{n}+
\frac{a(r_k)}{k}+\frac{4c_d \psi_d(n) m_k}{m_kn+r_k}.
\ee
Now, we take first the limit $k=m_kn+r_k$ to infinity
while $n$ is fixed. We obtain
\be{var-end}
\frac{a(n)}{n}\le \lim\inf \frac{a(k)}{k}\le
\lim\sup \frac{a(k)}{k}\le \frac{a(n)}{n}+
\frac{4c_d \psi_d(n) }{n}.
\ee
Then, we take $n$ to infinity to obtain the existence of
a limit for $a(k)/k$, say $\kappa(q,d)$. Looking at \reff{var-end}
with an identification of the limit, we have, for any $n$
\[
E[\L_q(n)]\le n\kappa(q,d)\le E[\L_q(n)]+4c_d \psi_d(n).
\]
and this is \reff{moy-key}.
\qed
\subsection{Variance Estimates}
We estimate now the variance of $\L_q(n)$, and prove \reff{key-var1}
and \reff{key-var2} of Theorem~\ref{theo-clt}. 

\noindent{\bf Step 1:} We show first that \reff{key-var1} holds in any
dimension greater or equal to 3.
To estimate the variance of $\L_q(n)$, we use the following simple fact.
If $X,Y,Z$ are random variables, and $\epsilon>0$, then
\be{var.1}
Y\le X\le Y+Z\Longrightarrow 
\var(X)\le (1+\epsilon) \var(Y)+(1+\frac{1}{\epsilon}) E[Z^2].
\ee
Indeed, we have $|X-E[Y]|\le |Y-E[Y]|+Z$ (note that $Z\ge 0$), and
\[
\var(X)=\inf_c\ E[(X-c)^2]\le E[(X-E[Y])^2]\le (1+\epsilon) E[(Y-E[Y])^2]+
(1+\frac{1}{\epsilon}) E[Z^2].
\]
Thus, we have from \reff{decomp.3} and \reff{var.1}
\be{var.2}
S_1\le \L_q(n)\le S_1+\I_1(n_1,n_2)\Longrightarrow
\var(\L_q(n))\le (1+\epsilon) \var(S_1)+(1+\frac{1}{\epsilon}) E[\I_1^2(n_1,n_2)]
\ee
Similarly as in \reff{var.3}, we have a constant $C_d$ such that
\be{var.3b}
E[\I_1^2(n_1,n_2)]\le C_d\psi_d^2(\max(n_1,n_2))\le C_d \psi_d^2(n),
\ee
where we only used that $\psi_d$ is increasing.
Thus, 
\be{var.4}
\var(\L_q(n))\le (1+\epsilon) \pare{ \var(\L_q(n_1))
+\var(\L_q(n_2))}+
(1+\frac{1}{\epsilon}) C_d\psi_d^2(n).
\ee
Now, when we choose the almost dyadic decomposition of Section~\ref{sec-d3},
\reff{quasi-dyadic}, and using induction, we have
\be{var.5}
\begin{split}
\var(\L_q(n))\le&(1+\epsilon)^L\pare{ \sum_{k=1}^{2^L}
\var(\L_q(n_k^{(L)}))}\\
&\qquad+(1+\frac{1}{\epsilon}) C_d'\sum_{k=1}^{2^L} (1+\epsilon)^{k-1}
2^{k-1} \psi_d^2(\frac{n}{2^{k-1}}).
\end{split}
\ee
Recall that $\psi_d^2(k)\le k$ for $d\ge 3$. Thus, when reaching $L=\lfloor\log_2(n)
\rfloor$,
$\var(\L_q(n_k^{(L)}))$ are of order 1, 
and there is a constant $C$, such that
\be{var.6}
\var(\L_q(n))\le C(1+\epsilon)^L 2^L+C_d'(1+\frac{1}{\epsilon})
\frac{(1+\epsilon)^L}{\epsilon} n.
\ee
Choosing $\epsilon=1/L$, we obtain \reff{key-var1} in $d\ge 3$. 

\noindent{\bf Step 2:} We consider now $d\ge 4$. We show that there is a constant
$C_d$ such that
\be{var.10}
\var(\L_q(n))\le C_d n.
\ee
We go back to \reff{var.4} and optimize over $\epsilon$ to obtain
\be{var.7}
\begin{split}
\var(\L_q(n))\le&\pare{ \var(\L_q(n_1))}+ \var(\L_q(n_2))+
C_d' \psi_d^2(\max(n_1,n_2))\\
&+2\pare{\pare{ \var(\L_q(n_1))
+\var(\L_q(n_2))}C_d'\psi_d^2(\max(n_1,n_2))}^{1/2}.
\end{split}
\ee
Now, choose first $n=2^m$, and $n_1=n_2=2^{m-1}$, and set $a_k=
\var(\L_q(2^k))2^{-k}$. Then, using \reff{key-var1} to estimate 
the cross-product in \reff{var.7}, we have
\be{var.8}
a_m\le a_{m-1}+r_m,\quad\text{with}\quad
r_m=\frac{C_d'\psi_d^2(2^m)}{2^m}+2\pare{\frac{C_d'c_d m^2 \psi_d^2(2^m)}{2^m}}^{1/2}
\ee
When $d\ge 4$, $\psi_d^2(2^m)\le C m^2$, and $\{r_m,m\in \N\}$ defines
a convergent series. Thus,
\be{var.9}
a_m\le a_0+\sum_{k=1}^{m} r_k\le c_d:=a_0+\sum_{k=1}^{\infty} r_k\Longrightarrow
\var(\L_q(2^m))\le c_d 2^m.
\ee
Now, write any integer $n$ in terms of its binary decomposition 
 $n=2^{m_1}+\dots+2^{m_k}$, with
$0\le m_1<m_2<\dots<m_k$. We call now $n_1=2^{m_k}$, and $n_2=n-n_1$, and note
that $n_1\ge n_2$.
In $d\ge 4$, we use the bound $\psi_d(k)\le \log(k)$
in \reff{var.7}, and the estimate \reff{key-var1} in bounding the
term $\var(\L_q(n_1))+\var(\L_q(n_2))$ which appears in the square root
in \reff{var.7}. Thus, we obtain that there exists a constant $c$ independent of $n$
such that
\be{var.26}
\var(\L_q(n))\le \var(\L_q(n_1))+\var(\L_q(n_2))+c m_k^2 \sqrt{ 2^{m_k}}.
\ee 
By iterating \reff{var.26}, we obtain using \reff{var.9}
\be{var.25}
\begin{split}
\var(\L_q(n))\le&\sum_{j=1}^k \var(\L_q(2^{m_j}))+c\sum_{j=1}^k m_j^2 \sqrt{ 2^{m_j}}\\
\le &c_d\sum_{j=1}^k  2^{m_j}+c\sum_{j=1}^k \frac{ m_j^2}{\sqrt{ 2^{m_j}}} 2^{m_j}\\
\le & (c_d+cc_3) n,
\end{split}
\ee
where $c_3$ is a constant such that for any $m$, $m\le c_3 \sqrt{ 2^{m}}$.

\noindent{\bf Step 3:} We show now how to obtain \reff{key-var2}. 
Note first that using similar arguments as those leading 
to \reff{var.2} and \reff{var.7}, we have
\be{var.11}
\pare{ \var(\L_q(n_1))+\var(\L_q(n_2))}\le \var(\L_q(n))+
C_d' \psi_d^2(\frac{n}{2})+2\pare{\var(\L_q(n))C_d' \psi_d^2(\frac{n}{2})}^{1/2}.
\ee
Thus, using \reff{key-var1} and \reff{var.11}, 
there is $c_1>0$ such that for any integer $j$, 
\be{var.12}
|\var(\L_q(2^{j}))- 2 \var(\L_q(2^{j-1}))|\le c_1 j\sqrt{ 2^j}.
\ee
Now, we consider $m,l,i$ integers, such that $2^m=2^l2^i$, and
consider for $j=1,\dots,l$ the system of inequalities obtained from \reff{var.12}
\be{var.13}
|2^j \var(\L_q(2^{i+l-j}))- 2^{j-1} \var(\L_q(2^{i+l-j+1}))|\le c_1(i+l-j+1)2^{j-1}
\sqrt{2^{i+l-j+1}}.
\ee
By summing \reff{var.13} for $j=1,\dots,l$, and using the triangle inequality, we obtain
\be{var.14}
|2^l \var(\L_q(2^{i}))-  \var(\L_q(2^{m}))|\le 
c_1\sqrt{2^{i+l}}\sum_{j=1}^l (i+l-j+1)\sqrt{2^{j-1}}.
\ee
By dividing both sides
of \reff{var.14} by $2^m$, we have a constant $c_2$ such that
\be{var.15}
|\frac{\var(\L_q(2^{i}))}{2^i}-\frac{\var(\L_q(2^{m}))}{2^m}|\le
\frac{c_2i\sqrt{2^l}} {\sqrt{2^{i+l}}}.
\ee
In \reff{var.15}, we take first the limit $l$ to infinity (recall that $2^m=2^l2^i$),
 then $i$ to infinity to conclude that there exists 
\be{var.16}
\lim_{n\to\infty} \var(\L_q(2^n))/2^n=v(q,d),
\quad\text{and}\quad |\frac{\var(\L_q(2^n))}{2^n}- v(q,d)|\le \frac{c_2 n}{\sqrt{2^n}}.
\ee
It is easy to conclude \reff{key-var2}. Indeed, for any integer $n$, consider its
dyadic decomposition, say $n=2^{m_1}+\dots+2^{m_k} $, and note that using
\reff{var.11} and Step 2, we can improve \reff{var.24} into
\be{var.17}
|\var(\L_q(n)) -\sum_{j=1}^k \var(\L_q(2^{m_j}))|\le c_1
\sum_{j=1}^k m_j\sqrt{2^{m_j}},
\ee
and \reff{var.16} allows us to conclude.

\subsection{The central limit theorem}\label{sec-lindeberg}
The aim of this section is to prove \reff{key-clt}.
We use the notations of Section~\ref{sec-d3}.
We fix $\delta_1>0$ small, and let $L_n$ be
the integer part of $\log_2(\sqrt{n} n^{-\delta_1})$. Note that this
choice $2^{L_n}\sim \sqrt{n} /n^{\delta_1}$ is different from the
choice of Section~\ref{sec-self} where $2^{L}\sim n^{1-\delta_0}$ 
for $\delta_0$ smaller that $2/(dq)$.

If we define $R(n)=\L_q(n)-S_q^{(L_n)}$, then \reff{low-eq.8} yields
\be{clt.eq3}
0\le R(n)\le \sum_{l=1}^{2^{L_n}} \I_l.
\ee
By subtracting to $\L_q(n)$ its average, we obtain
\be{clt.eq2}
\L_q(n)-E[\L_q(n)] =\sum_{k=1}^{2^{L_n}} Z^{(L_n)}_k+R(n)-E[R(n)],
\ee
with $Z^{(L_n)}_k=\L_q^{(k)}(n_k^{(L_n)})-E[\L_q^{(k)}(n_k^{(L_n)})]$.
As a first step, we show  that $R(n)/\sqrt{n}$ vanishes in law. More 
precisely, we show that
\be{clt.eq4}
\lim_{n\to\infty} \frac{E[R(n)]}{\sqrt{n}}=0.
\ee
Then, as a second step, we invoke the CLT for triangular arrays
(see for instance \cite{billingsley} p. 310), 
since we deal with independent random 
variables $\{Z^{(L_n)}_k,\ k=1,\dots,2^{L_n}\}$. The CLT 
states that for a standard normal variable $Z$
\be{CLT-array}
\frac{\sum_{k=1}^{2^{L_n}} Z^{(L_n)}_k}{\sqrt{\sum_{k=1}^{2^{L_n}} \var(Z^{(L_n)}_k)}}
 \stackrel{\text{law}}{\longrightarrow} Z,
\ee
provided that Lindeberg's condition holds. 
This latter condition reads in our context
\be{lindeberg}
\lim_{n\to\infty} \sup_{k\le 2^{L_n}} \frac{ 
E[\ind_{\acc{|Z^{(L_n)}_k|>\epsilon \sqrt{n}}}(Z^{(L_n)}_k)^2]}
{E[(Z^{(L_n)}_k)^2]}=0.
\ee
Assuming \reff{clt.eq4} and \reff{lindeberg} hold, we rely on Lemma~\ref{lem-moy}
to replace $E[\L_q(n)]$ by $n\kappa(q,d)$ at a negligible cost, and rely on
Theorem~\ref{theo-clt}
to replace the $\sum_k \var(Z^{(L_n)}_k)$ by $n v(q,d)$. 
Indeed, note that by \reff{key-var2}
\be{approx-var1}
|\var(Z^{(L_n)}_k)-n_k^{(L_n)} v(q,d)|\le c(q,d) \log(n_k^{(L_n)}) 
\sqrt{n_k^{(L_n)}},
\ee
so that by summing over $k=1,\dots,,2^{L_n}$,
\be{approx-var2}
|\sum_{k=1}^{2^{L_n}} \var(Z^{(L_n)}_k) -n v(q,d)|\le c(q,d) 2^{L_n}\sqrt{
\frac{n}{2^{L_n}}} \log\pare{\frac{n}{2^{L_n}}}\le c(q,d) n^{\frac{3}{4}}\ \ \frac{
\log( \sqrt{n} n^{\delta_1})}{\sqrt{n^{\delta_1}}}.
\ee

\noindent{\bf Step 1:} 
We estimate the expectation of $R(n)$. From \reff{decomp.5} and Lemma~\ref{lem-app.1},
with $b_i=i$,
\be{clt.6}
E[\I_l]\le \sum_{k=1}^{2^{l}}\sum_{i\ge 0} 2^q(i+1)^{q-1}e^{-\kappa_d i}
C_d\psi_d(n^{(l)}_k)\le C'_d 2^l \log\pare{\frac{n}{2^l}}.
\ee
Thus, $E[R(n)]\le C' 2^{L_n} \log(n)\le C'\frac{\log(n) \sqrt{n}}{n^{\delta_1}}$, and
$\lim_{n\to\infty} E[\frac{R(n)}{\sqrt{n}}]=0$.

\noindent{\bf Step 2:} To check Lindeberg's condition, we start with
estimating $P(|Z^{(L_n)}_k|\ge \epsilon \sqrt{n})$. To simplify
notation, we set $n_k=n^{(L_n)}_k$, and we note that
\be{clt-quest}
P(|Z^{(L_n)}_k|\ge \epsilon \sqrt{n})=P\pare{|\L_q(n_k)-
E[\L_q(n_k)]|\ge \xi_{n_k} n_k},\ \text{and,}\
\xi_{n_k} =\frac{\epsilon \sqrt{n}}{n_k}
\ge \frac{\epsilon}{2n_k^{\delta}},
\ee
with $\delta=\frac{2\delta_1}{1+2\delta_1}$. Thus, Lindeberg's condition
is written as a large deviation for $\L_q(n_k)$.
Note that $n_k$ is almost the scale of the CLT.
We now use Remark~\ref{rem-xin}, and Lemma~\ref{lem-app4} of the Appendix.
We apply \reff{eq-xin1}, \reff{eq-xin2} of Remark~\ref{rem-xin}, and
\reff{clt.7} and \reff{clt.8} of Lemma~\ref{lem-app4}, to obtain
for arbitrarily small $\epsilon'$ and $\delta$
\be{clt.10}
\begin{split}
P\pare{|Z^{(L_n)}_k|\ge \epsilon \sqrt{n}}&\le 
P\pare{Z^{(L_n)}_k\ge \epsilon \sqrt{n}}+
P\pare{Z^{(L_n)}_k\le -\epsilon \sqrt{n}}\\
&\le \exp\pare{ -C\pare{\frac{2\epsilon}{
n_k^\delta}}^{\max(\frac{1}{q},\frac{2}{d}\gamma)+\frac{2}{d}}
n_k^{\min(1/q_c(d),1/q)-\epsilon'} }+
e^{-\frac{\epsilon}{4}n_k^{1-q\delta_0-\delta}}.
\end{split}
\ee
Inequality \reff{clt.10} with the
uniform bound $|Z^{(L_n)}_k|\le n^{q(\delta+\frac{1}{2})}$,
and the lower bound on $\var(Z^{(L_n)}_k)$ in \reff{approx-var1},
imply that Lindeberg's condition \reff{lindeberg} holds.

\section{Appendix}\label{sec-append}
In this section, we recall and improve some key estimates for
dealing with large deviation for intersection local times.
First, we recall a special form of Lemma 5.1 of \cite{A06}.
\bl{lem-5.1}[Lemma 5.1 of \cite{A06}]
Assume $\{Y_1,\dots,Y_n\}$ are positive and independent. Furthermore,
assume that there is a constant $C>0$ such that 
for any $i\in \{1,\dots,n\}$
\be{normalize-geo}
\forall t>0\quad P(Y_i>t)\le C\exp(-t).
\ee
Then, for some $c_u>0$, and any $0<\delta<1$, we have for any integer $n$
\be{prel.6}
P\pare{ \sum_{i=1}^n  \pare{Y_i-E[Y_i]} \ge x_n}\le
\exp\pare{ c_u \delta^{2(1-\delta)} 
n\max_{i}\pare{E[Y_i^2],E[Y_i^2]^{1-\delta}}-\frac{\delta}{2} x_n}.
\ee
\el
Secondly, we improve Lemma 5.3 of \cite{A06} into inequalities we
believe are optimal.
Consider two independent random walks $\{S(n),\tilde S(n),\ n\in \N\}$, 
and for an integer $k$, 
denote $\tilde D_n(k):=\{z\in \Z^d:\ \tilde l_n(z)>k\}$. We recall
that if $l_n$ is the local times and $A$ a subset of $\Z^d$, 
then $l_n(A)=\sum_{z\in A} l_n(z)$.
\bl{lem-app.1}
Assume dimension $d\ge 3$. There are constants 
$C_d,C_d',\kappa_d$ such that
\be{app-key1}
E\cro{ l_{n}(\tilde D_n(k)) }\le C_d e^{-\kappa_d k}\psi_d(n),
\quad\text{with}\quad \psi_d(n)=
\left\{ \begin{array}{ll}
n^{1/2} & \mbox{ for } d=3 \, , \\
\log(n) & \mbox{ for } d=4 \, , \\
1 & \mbox{ for } d\ge 5 \, ,
\end{array} \right.
\ee
and,
\be{app-key2}
E\cro{ l_{n}(\tilde D_n(k))^2 }\le C_d' e^{-\kappa_d k}\psi_d(n)^2.
\ee
\el
Finally, we prove the following lemma. This result is
not optimal, but suffices for our purpose. 
\bl{lem-app4} Assume $d\ge 3$, and take
$1>\xi_n\ge n^{-\delta}$ for $\delta\le \delta_0/3$ small enough.
(i) when $q\ge q_c(d)$, then for any $\epsilon>0$,
\be{clt.7}
P\pare{ \L_q(n)-E[\L_q(n)]\ge \xi_n n}\le \exp\pare{
-C \xi_n^{\frac{1}{q}+\frac{2}{d}} n^{\frac{1}{q}-\epsilon}}.
\ee
(iii) For any $q>1$,
\be{clt.8}
P\pare{ \L_q(n)-E[\L_q(n)]\le -\xi_n n}\le \exp\pare{ -\frac{\xi_n}{2}
n^{1-q\delta_0}}.
\ee
\el
\subsection{Proof of Lemma~\ref{lem-app.1}}
To emphazise the starting point, we denote by $P_z$ the law
of the random walk started at site $z\in \Z^d$.
We let $H_z=\inf\{n\ge 0:\ S(n)=z\}$, and
use Theorem 3.2.3 of Lawler~\cite{lawler}. 
\be{app.1}
\sum_{z\in \Z^d} P_0(H_z\le n)^2\le \sum_{z\in \Z^d} \pare{ \sum_{k=0}^n
P_0(S(k)=z)}^2\le C_d \psi_d(n).
\ee
Now call $P_0(l_\infty(0)>1)=e^{-\kappa_d}<1$, the return probability,
and 
\[
E_0[l_{\infty} (0)]=\frac{1}{1-e^{-\kappa_d}},\quad\text{and}\quad
E_0[l_{\infty} (0)^2]=\frac{1+e^{-\kappa_d}}{(1-e^{-\kappa_d})^2}
\]
It is easy to see that for any $z\in \Z^d$
\[
E_0\cro{l_n(z)}\le P_0(H_z\le n) E_0[l_{\infty} (0)]
,\quad\text{and}\quad E_0\cro{l_n^2(z)}\le P_0(H_z\le n) E_0[l_{\infty}^2(0)].
\]
Similarly,
\[
P_0\pare{ l_n(z)>k} \le P_0\pare{ H_z\le n}P_z\pare{ l_\infty(z)>k}=
e^{-\kappa_d k} P_0\pare{ H_z\le n}.
\]
Thus, there is $C_d$ such that
\ba{app.2}
E\cro{ l_n(\tilde D_n(k))}&=&\sum_{z\in \Z^d} E_0\cro{ l_n(z)}
P_0\pare{ l_n(z)>k}\cr
&\le&e^{-\kappa_d k} E_0[l_{\infty} (0)] \sum_{z\in \Z^d}  P_0\pare{ H_z\le n}^2
\le C_d  e^{-\kappa_d k} \psi_d(n).
\ea
Now, we expand the square of $l_n(\tilde D_n(k))$
\ba{app.3}
 l_n(\tilde D_n(k))^2&=&\pare{ \sum_{z\in \Z^d}  l_n(z)\ind\acc{\tilde 
l_n(z)>k}}^2\cr
&=& \sum_z l_n(z)^2 \ind\acc{\tilde l_n(z)>k}+
\sum_{z\not= z'} l_n(z)l_n(z')\ind\acc{\tilde l_n(z)>k, \tilde l_n(z')>k}.
\ea
After taking the expectation of $l_n(\tilde D_n(k))^2$
\ba{app.4}
E\cro{ l_n(\tilde D_n(k))^2}&=&\sum_z E_0\cro{ l_n(z)^2 } P_0\pare{
l_n(z)>k}+\sum_{z\not= z'}  E_0\cro{ l_n(z)l_n(z')}
P_0\pare{ l_n(z)\wedge l_n(z')>k}\cr
&\le & E_0\cro{ l_n(0)^2 } e^{-\kappa_d k} \sum_z P_0\pare{ H_z\le n}^2\cr
&&\quad +
\sum_{z\not= z'}  E_0\cro{ l_n(z)l_n(z')}P_0\pare{ l_n(z)\wedge l_n(z')>k}.
\ea
Now, in the last term in \reff{app.4},
we distinguish which of $z$ or $z'$ is hit first.
\be{app.5}
\begin{split}
P_0\pare{ l_n(z)\wedge l_n(z')>k}\le& P_0\pare{ H_z<H_{z'},\ l_n(z')>k}
+P_0\pare{ H_{z'}<H_{z},\ l_n(z)>k}\\
\le & P_0\pare{ H_z\le n}P_z\pare{ l_n(z')>k}+
P_0\pare{ H_{z'}\le n}P_{z'}\pare{ l_n(z)>k}\\
\le & e^{-\kappa_d k} \pare{ P_0\pare{ H_z\le n}P_z\pare{ H_{z'}\le n}+
P_0\pare{ H_{z'}\le n}P_{z'}\pare{ H_{z}\le n}}.
\end{split}
\ee
We treat now the term $E_0\cro{ l_n(z)l_n(z')}$. We have
\be{app.6}
\begin{split}
E_0\cro{ l_n(z)l_n(z')}=& \sum_{k<k'\le n} \pare{ P_0(S(k)=z)P_z(S(k'-k)=z')+
P_0(S(k)=z')P_{z'}(S(k'-k)=z)}\\
\le & E_0\cro{l_n(z)}E_{z}\cro{l_n(z')}+E_0\cro{l_n(z')}E_{z'}\cro{l_n(z)}\\
\le & E_0[l_{\infty} (0)]^2\pare{ P_0\pare{ H_z\le n}P_z\pare{ H_{z'}\le n}+
P_0\pare{ H_{z'}\le n}P_{z'}\pare{ H_{z}\le n}}.
\end{split}
\ee
Thus, with the help of \reff{app.5} and \reff{app.6}, \reff{app.4} reads
\be{app.7}
\begin{split}
E\cro{ l_n(\tilde D_n(k))^2}\le & E_0\cro{ l_n(0)^2 } 
e^{-\kappa_d k} \sum_z P_0\pare{ H_z\le n}^2\\
+ E_0[l_{\infty} (0)]^2 e^{-\kappa_d k} &\sum_{z\not= z'} 
\pare{ P_0\pare{ H_z\le n}P_z\pare{ H_{z'}\le n}+
P_0\pare{ H_{z'}\le n}P_{z'}\pare{ H_{z}\le n}}^2\\
\le &
E_0\cro{ l_n(0)^2 } e^{-\kappa_d k} \sum_z P_0\pare{ H_z\le n}^2\\
+2E_0[l_{\infty} (0)]^2 e^{-\kappa_d k}&\sum_{z\not= z'} 
P_0\pare{ H_z\le n}^2P_z\pare{ H_{z'}\le n}^2+
P_0\pare{ H_{z'}\le n}^2P_{z'}\pare{ H_{z}\le n}^2.
\end{split}
\ee
Now, we use translation invariance and \reff{app.1}
\[
\sum_{z\not= z'} P_0\pare{ H_z\le n}^2P_z\pare{ H_{z'}\le n}^2\le
\pare{\sum_z P_0\pare{ H_z\le n}^2}^2\le C_d^2\psi_d(n)^2.
\]
The result \reff{app-key2} follows at once.

\subsection{Proof of Lemma~\ref{lem-app4}}
The proof of (i) follows from \reff{super-crit} of
Lemma~\ref{lem-topI}, and Remark~\ref{rem-superxi-n} which deals
with the contribution of $\{z:\ l_n(z)<\xi_n^{1/q} n^{1/q-\epsilon}\}$.
Using that for a transient walk, the local time of a site is
bounded by a geometric variable, we have for a small $\delta>0$ and
a constant $c>0$,
\[
P\pare{ \sum_{z} \ind\acc{l_n(z)
\ge \xi_n^{1/q}n^{1/q-\epsilon}}\ l_n^q(z)\ge n\xi_n\delta}\le
P\pare{\exists z: l_n(z) \ge \xi_n^{1/q}n^{1/q-\epsilon}}\le
ne^{-c\xi_n^{1/q}n^{1/q-\epsilon}}.
\]

Point (ii) follows from  the lower bound
in \reff{low-eq.8}: $\L_q(n)\ge S_q^{(L)}$. 
Indeed, we choose $\delta=\delta_0/3$, 
(with $\delta_0<2/(dq)$) and $L$ such that
$2^L\sim n^{1-\delta_0}$. Then, we first have
\[
\L_q(n)-E[\L_q(n)]\le -\xi_n n\Longrightarrow S_q^{(L)}-
E[\L_q(n)]\le -\xi_n n.
\]
Now, Lemma~\ref{lem-LD} gives us
\be{clt.9}
P\pare{ S_q^{(L)}-E[\L_q(n)]\ge -\xi_n n}
\le \exp\pare{ -\frac{\xi_n}{2}n^{1-q\delta_0}}.
\ee

\noindent{\bf Acknowledgements}. We would like to thank
an anonymous referee for his questions, suggestions and comments which
led to correcting a mistake in the first version, as well
as an improved exposition.

\end{document}